\documentclass[11pt]{article}

\usepackage{array} 
\usepackage{amsmath} 
\usepackage{amssymb} 
\usepackage{amsthm} 
\usepackage{cases} 
\usepackage{diagbox} 
\usepackage{enumitem}
\usepackage{extarrows} 
\usepackage{framed} 
\usepackage{indentfirst}
\usepackage{geometry} 
\usepackage{graphicx} 
\usepackage{makecell} 
\usepackage{mathrsfs} 
\usepackage{mathtools} 
\usepackage{pdfpages}
\usepackage{ulem} 
\usepackage{upgreek} 
\usepackage{verbatim}
\usepackage{wrapfig}
\usepackage{xcolor,color} 
\usepackage[framemethod=TikZ]{mdframed}

\numberwithin{equation}{section}

\theoremstyle{plain}
\newtheorem{theorem}{Theorem}[section]

\newtheorem{proposition}[theorem]{Propsition}
\newtheorem{corollary}[theorem]{Corollary}

\theoremstyle{definition}

\theoremstyle{remark}

\newtheorem*{prof}{\textit{Proof}}

\def\<{\langle}\def\>{\rangle}

\allowdisplaybreaks[4]

\newcommand*\Abstract[1]{%
	\begingroup\noindent\leftskip=.7cm
	\rightskip\leftskip
	\textbf{Abstract}\quad #1%
	\par\vspace*{1mm}\endgroup}
\newcommand*\Keywords[1]{%
	\begingroup\noindent\leftskip=.7cm
	\rightskip\leftskip
	{\hangafter=1\hangindent=2.1cm\noindent\textbf{Key words}\quad #1%
		\par\vspace*{1mm}}\endgroup}
\newcommand*\MRSubClass[1]{%
	\begingroup\noindent\leftskip=.7cm
	\rightskip\leftskip
	\textbf{2010 MR Subject Classification}\quad #1%
	\par\vspace*{1mm}\endgroup}

\begin{document}
	
	\bigskip\bigskip
	\noindent{\Large\bf Limit theorems for supercritical remaining-lifetime age-structured branching processes}
	
	\vspace*{1mm}
	
	\noindent{
		
		Ziling Cheng
	}
	\vspace*{1mm}
	
	\Abstract{We study supercritical age-structured branching models starting from a single particle with a random lifetime, where the reproduction law depends on the remaining lifetime of the parent.  The lifespan of an individual is decided at its birth and its remaining lifetime decreases at the unit speed. A necessary and sufficient condition is provided for the convergence of the Malthusian normalized random measures. The Malthusian type limit theory in a functional form can be strengthened to hold with probability one under some ``$L\log L$'' conditions. We further prove a central limit theory with a random normalization factor.}
	
	
	\MRSubClass{60J80, 60J85, 60F15, 60F05} 
	
	\section{Introduction}
	
	Limit theorems form an important and fundamental field of the theory of supercritical branching processes. The longtime behavior of the population size of different branching models has been investigated by several authors. For the supercritical Galton-Watson process $(Z_n  : n\ge 0)$, 
	Seneta \cite{Seneta68} proved that there exists a series of positive constants $\{c_n; n\ge 0\}$ such that the scaled population sizes $c_n Z_n$ have a non-degenerate limit distribution. Heyde \cite{Heyde70} then showed that the existence can be strengthened to hold almost surely. In particular, let $m$ be the expected number of offspring per particle and $c_n=m^{-n}$. Using the theory of positive martingale, Kesten and Stigum \cite{Kesten66} gave a necessary and sufficient condition (called the \textit{Kesten-Stigum $L \log L$ criterion}) for the almost sure existence of a non-degenerate limit distribution. Based on this, Athreya \cite{Athreya68} considered the multitype continuous time Markov branching processes $(Z_t : t\ge 0)$. They established the almost sure convergence results for the normalized population sizes $e^{-\lambda_1 t} Z_t$, where $\lambda_1$ is the maximal real eigenvalue of the infinitesimal generator  of the mean matrix semigroup. Asmussen and Hering \cite{Asmussen76} further generalized these results to the case of branching Markov processes under some $L \log L$ conditions. In the case of the age-dependent branching process (or simply the \textit{B-H process}) $(Z_t:t\ge 0)$ introduced by Bellman and Harris \cite{Bellman52}, Cohn \cite{Cohn82} showed that there always exists constants $C_{\mathrm{t}}$ such that $Z_t / C_t$ almost surely converge to some non-degenerate random variable. Schuh \cite{Schuh82} identified $C_t$ to be as in \cite{Seneta68} for the Galton-Watson process. The necessary and sufficient condition ($L \log L$ criterion) for the convergence in law of $e^{-\alpha t}Z_t $ was given in Athreya \cite{Athreya69}, where $\alpha$ is the Malthusian parameter of the B-H process. Using this result and the limit behavior of the age distribution, Athreya \cite{Athreya76} showed that under $L \log L$ criterion  $e^{-\alpha t}Z_t$ converge almost surely to a non-degenerate limit, which extending the Kesten-Stigum theorem to the age-dependent case. Similar results have been established for several supercritical branching processes; see also \cite{Biggins04, Conner66, Conner67, Harris63, Harris59, Hering73, Kharlamov68, Kurtz74, Kurtz97, Lyons95}.
	
	The age-dependent birth and death branching process (or simply the \textit{C-M-J process}) $(Z_t : t\ge 0)$ introduced by Crump and Mode \cite{Crump68} and Jagers \cite{Jagers69} is a more general branching model. Several authors have also studied extending the Kesten-Stigum theorem to the general class of branching processes. Crump and Mode \cite{Crump68, Crump69} discussed the convergence in mean square of $Z_t/\mathbf{E}Z_t$. Doney \cite{Doney72a, Doney72b} found weaker conditions ($L \log L$ criterion) for the convergence in distribution of $Z_t/\mathbf{E}Z_t$. Benoît \cite[Chapter 5]{Benoît16} gave a new proof of the almost sure convergence of the normalized population size and then established a central limit theorem of the population size. Furthermore, Doney \cite{Doney76} extended the results in \cite{Doney72b} to the multi-type general age-dependent branching processes. In fact, in addition to the population size, the extreme behavior of age structure also has research significance. Individuals in the C-M-J process can be counted by the values of a random characteristic as proposed by Jagers \cite{Jagers74, Jagers75}. The limit theory for the random characteristic model in the supercritical case was also developed by Jagers \cite{Jagers74} under some conditions of second moment. Nerman \cite{Nerman81} established the convergence in probability of the Malthusian normalized supercritical C-M-J processes counted with a random characteristic under some mild regularity conditions. They further proved some almost sure convergence results, provided the tail of the reproduction point process and the characteristic both satisfy mild regularity conditions. The case where random characteristics depend not only on age but also on absolute time was considered by Jagers and Nerman \cite{Jagers84}. \textcolor{black}{Jagers \cite{Jagers84'} further considered the models with more general characteristics.}
	
	In a previous paper \cite{Cheng23+} we studied a class of remaining-lifetime age-structured branching processes with reproduction law depending on the remaining lifetime of the parent. The lifespan of an individual is distributed arbitrarily and determined at its birth and its remaining lifetime decreases at the unit speed. The model can be seen as a special case of the C-M-J process. In this paper, we consider a special case of the model, where the population descends from a single individual with a random lifetime and the reproduction regime is supercritical. It is obvious that in such a process the generation sizes form a Galton-Watson process, i.e., the so-called \textit{embedding Galton-Watson process}. 
	Some preliminary results are stated in Section 2, including some characterizations of the embedding Galton-Watson process. Let $\tilde{\alpha}$ be the strictly positive (and finite) Malthusian parameter. 
	
	We first give a sufficient and necessary condition for  the convergence of the Malthusian normalized random measures in Section 3 and 4. Then in Section 5, we strengthen the Malthusian type limit theory in a functional form to hold almost surely under some ``$L\log L$'' conditions
	. \textcolor{black}{Actually, we use a different method from Nerman \cite{Nerman81} and Jagers \cite{Jagers84'} to prove the almost surely convergence. We give a ``$L\log L$'' condition by considering the ``$L\log L$''-moment of the process $(X_t :t\ge 0)$.}
	
	We further establish a central limit theory with a random
	normalization factor (naturally conditioning on the non-extinction event) in Section 6. More precisely, we show that $(\langle X_t,f\rangle-\mathbf{E}[\langle X_t,f\rangle])/\sqrt{X_t (\infty)}$ converges in distribution to a normal random variable under additional (second moment) assumptions, where random variable $X_t (\infty)$ denotes the total number of particles alive at time $t\ge 0$. \textcolor{black}{In particular, Iksanov \textit{et al.} \cite{Iksanov21} considered a supercritical Crump-Mode-Jagers process $(\mathcal{Z}_t^{\varphi}: t \geq 0)$ counted with a random characteristic $\varphi$. They also proved a central limit theorem for $(\mathcal{Z}_t^{\varphi}: t \geq 0)$. More precisely, they showed that $(\mathcal{Z}_t^{\varphi}-a e^{\alpha t} W-H(t)) / \sqrt{t^k e^{\alpha t}}$ converges in distribution to a normal random variable with random variance, where $aW$ is the almost sure limit of $e^{-\alpha t}\mathcal{Z}_t^{\varphi}$, $k$ is a constant and $H(t)$ is a function. But the central limit theorem proved by us not be directly obtained from the above result. Indeed, we could not give a characterization of the joint distribution of random vector 
		$$\Big(\frac{\langle X_t,f\rangle-\mathbf{E}[\langle X_t,f\rangle]}{\sqrt{e^{\tilde{\alpha} t}}},\widetilde{W}_{\infty}^{1}\Big),$$ 
		where $\widetilde{W}_{\infty}^{1}$ is the almost sure limit of $e^{-\tilde{\alpha} t}X_t (\infty)$. This is due to the dependence of the two random variables.}
	
	Let $\mathbb{N}=\{0,1,2,\ldots\}$. Let $\mathfrak{M}(0,\infty)$ be the set of finite Borel measures on $(0,\infty)$ with the weak convergence topology. 
	Let $\mathfrak{N}(0,\infty)$ be the subset of $\mathfrak{M}(0,\infty)$ consisting of integer-valued measures. Let $\mathcal{B}(0,\infty)$ denote the Borel $\sigma$-algebra on $(0,\infty)$. Let $B(0,\infty)$ be the Banach space of bounded Borel functions on $(0,\infty)$ furnished with the supremum norm $\|\cdot\|$. Let $C(0,\infty)$ be the set of continuous functions in $B(0,\infty)$. Let $C^{1}(0,\infty)$ be the set of functions in $C(0,\infty)$ with bounded continuous derivatives of the first order. We use the superscript ``+'' to denote the subsets of positive elements
	, e.g., $B(0,\infty)^{+}$, 
	etc. For any function $f$ on $A\subset\mathbb{R}$, we understand that $f(x)=0$ for $x\in\mathbb{R}\backslash A$ by convention. For any $f\in B(0,\infty)$ and $\nu\in\mathfrak{M}(0,\infty)$ write $\langle\nu,f\rangle= \int_{(0,\infty)}f(x)\nu(d x)$. In the integrals, we make the convention that, for $a\le b\in\mathbb{R}$,
	\begin{align*}
		\int_{a}^{b}=\int_{(a, b]}\quad \text { and } \quad \int_{a}^{\infty}=\int_{(a, \infty)}.
	\end{align*}
	Let $\xrightarrow{d}$, $\xrightarrow{\mathbf{P}}$ and $\xrightarrow{\mathbf{P}\text{-a.s.}}$ stand for convergence in distribution, in probability with respect to $\mathbf{P}$ and almost surely with respect to $\mathbf{P}$, respectively.
	
	\section{An age-structured branching process}\label{section2}
	
	Let $\alpha\in C^{1}(0,\infty)^{+}$. For each $x\in (0,\infty)$, let $\{p(x, i):i\in\mathbb{N}\}$ be a discrete probability distribution with generating function
	\begin{align*}
		g(x,z)=\sum\limits_{i=0}^{\infty}p(x,i) z^{i}, \quad z \in[0,1].
	\end{align*}
	Then for any $x\in(0,\infty)$ let $g'(x,z)=\frac{\partial}{\partial z} g(x,z)$ and $g''(x,z)=\frac{\partial^2}{\partial z^2} g(x,z)$ denote respectively the first and second derivative with respect to $z$ of $g(x,z)$. Throughout the paper, we assume that $p(\cdot,i) \in C^{1}(0,\infty)^{+}$ for every $i \in \mathbb{N}$ and
	\begin{align}\label{2.1}
		\|g'(\cdot, 1-)\|=\sup _{x>0} \sum_{i=1}^{\infty} p(x, i) i<\infty,
	\end{align}
	Then we have $g(\cdot, z)\in C^{1}(0,\infty)^{+}$ for each $z\in [0,1]$. Let $\beta=\|\alpha g'(\cdot,1-)\|<\infty$. Let $G$ be a continuous probability distribution on $(0,\infty)$. 
	A branching particle system is characterized by the following properties:
	
	\begin{description}
		
		\item[\textmd{(i)}] It starts with a single particle born at time $t=0$ with lifetime $L$, where $L$ is a random variable taking values in $(0,\infty)$ with distribution $G$.
		
		\item[\textmd{(ii)}] The remaining lifetimes of the particles decrease at the unit speed, i.e., they move according to realizations of the deterministic process $\xi=(\xi_{t}\vee 0)_{t \ge 0}$ in $(0,\infty)$ defined by $\xi_{t}=\xi_{0}-t$.
		
		\item[\textmd{(iii)}] A particle gives birth to offspring during its life. For a particle which is alive at time $r \ge 0$ with remaining lifetime $x>t-r>0$, the conditional probability of having not given birth by time $t$ is $\exp\{-\int_{0}^{t-r} \alpha(x-s) d s\}$.
		
		\item[\textmd{(iv)}] When a particle gives birth at remaining lifetime $x > 0$, it firstly gives birth to a random number of offspring according to the probability law $\{p(x, i): i\in \mathbb{N}\}$ determined by the generating function $g(x, \cdot)$, those offspring then choose their life-lengths in $(0,\infty)$ independently of each other according to the continuous probability distribution $G(d t)$.
		
	\end{description}
	
	Notice that (\ref{2.1}) ensures the number of offspring born at a branching event as above is almost surely finite.  We assume that the lifetimes and the offspring reproductions of different particles are independent. Let $X_{t}(B)$ denote the number of particles alive at time $t \ge 0$ with remaining lifetimes belonging to the Borel set $B \subset (0,\infty)$. Then $(X_{t}: t\ge 0)$ is a Markov process with state space $\mathfrak{N}(0,\infty)$. Suppose that the process is defined on a filtered probability space $(\Omega, \mathcal{F}, \mathcal{F}_{t}, \mathbf{P})$ satisfying the usual hypotheses. Let $\hat{X}=(\Omega, \mathcal{F}, \mathcal{F}_{t}, \hat{X}_{t}, \hat{\mathbf{P}}_{x})_{x\in (0,\infty)}$ be a c\'{a}dl\'{a}g realization of the $(\alpha,g,G)$-process starting from a single particle with lifetime $x\in(0,\infty)$. Such a process introduced in Cheng and Li \cite{Cheng23+} is characterized by the above properties (ii)-(iv). We further assume that $L$ and $\hat{X}$ are independent. Then it is easy to obtain
	\begin{align}\label{2.3}
		\mathbf{P}(X_t \in\cdot)=\int_{0}^{\infty} \hat{\mathbf{P}}_{x}(\hat{X}_t\in\cdot) G(dx),\quad t\ge 0.
	\end{align}
	
	We refer to Cheng and Li \cite{Cheng23+} for the formulation of $(\alpha,g,G)$-processes.	Notice that the remaining lifetime of a living particle must be greater than $0$, then the above properties imply that
	\begin{align}\label{2.4}
		\mathbf{E}[\exp \{-\langle X_{t}, f\rangle\}]=\langle G, e^{-u_{t} f}\rangle, \quad f \in B(0, \infty)^{+},
	\end{align}
	where $u_t f(x)$ is the unique solution to the following renewal equations:
	\begin{align}\label{2.5}
		e^{-u_{t} f(x)}=e^{-f(x-t)}+\int_{0}^{t}\big[g(x-s,\langle G, e^{-u_{t-s} f}\rangle)-1\big] e^{-u_{t-s} f(x-s)} \alpha(x-s) d s,
	\end{align}
	or
	\begin{align}\label{2.5'}
		u_{t}f(x)=f(x-t)+\int_{0}^{t} \alpha(x-s)\big[1-g(x-s,\langle G,e^{-u_{t-s} f}\rangle)\big]ds.
	\end{align}
	Using Cheng and Li \cite[Proposition 2.2]{Cheng23+} and dominated convergence theorem, we naturally obtain that	
	\begin{proposition}\label{prop2.1}
		For any $t \ge 0$ we have
		\begin{align}\label{2.1'}
			\mathbf{E}[\langle X_{t}, f\rangle]=\langle G, \pi_{t} f\rangle, \quad f \in B(0,\infty),
		\end{align}
		where $(\pi_{t})_{t \ge 0}$ is the semigroup of bounded kernels on $(0,\infty)$ defined by
		\begin{align}\label{2.2'}
			\pi_{t} f(x)=f(x-t)+\int_{0}^{t} \alpha(x-s) g'(x-s,1-) \langle G,\pi_{t-s}f\rangle ds.
		\end{align}
	\end{proposition}
	
	
	By integrating both sides of (\ref{2.2'}) with respect to $G(dx)$ we get
	\begin{align}\label{2.3'}
		\langle G,\pi_{t} f\rangle=\int_{0}^{t} \langle G,\pi_{t-s}f\rangle F(ds)+\int_t^{\infty}f(x-t) G(dx),
	\end{align}
	where
	\begin{align*}
		F(ds):=\Big(\int_{s}^{\infty}\alpha(x-s) g'(x-s,1-) G(dx)\Big)ds.
	\end{align*}
	By (\ref{2.3'}) and the general result on defective renewal equation; see, e.g., Jagers \cite[Theorem 5.2.8]{Jagers75}, we have the following proposition.
	\begin{proposition}\label{prop2.2} 
		Suppose that $m:=\int_{0}^{\infty}F(ds)\in (1,\infty)$. Then there exists a unique constant $\tilde{\alpha}\in (0,\infty)$ such that
		\begin{align}\label{2.*}
			\int_0^{\infty} e^{-\tilde{\alpha}t}F(dt)=1
		\end{align}
		and for any $f\in B(0,\infty)^{+}$,
		\begin{align}\label{2.*'}
			\lim_{t\rightarrow \infty} e^{-\tilde{\alpha} t} \langle G,\pi_t f \rangle=\frac{\int_0^{\infty} e^{-\tilde{\alpha} u} d u \int_u^{\infty} f(x-u) G(d x)}{\int_0^{\infty} u e^{-\tilde{\alpha} u} F(d u)}=:a(f) <\infty.
		\end{align}
		(Here we understand that $a(f) =0$ if $\int_0^{\infty} t e^{-\tilde{\alpha} t} F(d t)=+\infty$.)
	\end{proposition}
	Indeed, (\ref{2.*'}) holds since $e^{-\tilde{\alpha}t}\int_t^{\infty} f(x-t) G(d x)$ is directly Riemann integrable over $[0,\infty)$. In this paper we deal with the \textit{supercritical} branching case, which means that
	\begin{align}\label{super}
		1<m<\infty \quad \text{and} \quad \int_0^{\infty} t e^{-\tilde{\alpha} t} F(d t)<\infty.
	\end{align}
	Heuristically, let $N_{x}(t)$ be the number of children born to the ancestor with lifetime $x\in (0,\infty)$ in the time interval $(0,t]$. Let $M(d t, d u, d n, d v)$ be an $(\mathcal{F}_{t})$-Poisson random measure on $(0, \infty)^{2} \times \mathbb{N} \times(0,1]$  with intensity $d t d u \pi(dn) d v$, where  $\pi(d n)$ denotes the counting measure on $\mathbb{N}$. By arguments similar to Cheng and Li \cite[(3.1)]{Cheng23+} we have
	\begin{align}\label{4.1'}
		N_{x}(t)=\int_{0}^{x\wedge t}\int_{0}^{\alpha (x-s)}\int_{\mathbb{N}}\int_{0}^{p(x-s,n)} n\ M(ds, du, dn, dv),\quad t\ge 0.
	\end{align}
	It is easy to see that for any $n\in\mathbb{N}$ and $t\ge 0$ we have
	\begin{align*}
		\mathbf{P}(N_{L}(t)=n)=\int_{0}^{\infty} \mathbf{P}(N_{x}(t)=n)\ G(dx).
	\end{align*}
	Then we obtain
	\begin{align*}
		\mathbf{E}[N_{L}(t)]=\int_{0}^{\infty} G(dx) \int_{0}^{x\wedge t} \alpha(x-r) g'(x-s,1-) ds,
	\end{align*}
	which implies that
	\begin{align*}
		m=\mathbf{E}[N_{L}(\infty)]=\int_{0}^{\infty} G(dx) \int_{0}^{x} \alpha(x-s) g'(x-s,1-) ds,
	\end{align*}
	where $N_{L}(\infty)$ is the total number of children born to the ancestor of $(X_{t}: t\ge 0)$ during its life. For simplicity, we then  write $N(t):=N_{L}(t)$.
	
	We now discuss the extinction probability. To this end, we need to consider the generating function of $X$. For any $f\in B(0,\infty)^{+}$, $\langle X_t ,f\rangle$ can be represented as the sum
	\begin{align}\label{2.14}
		\langle X_t ,f\rangle=\sum_{j=1}^{N(t)} \langle X_{t-t_j}^{(j)} ,f\rangle + f(L-t),
	\end{align}
	involving $N(t)$ independent daughter processes $(X_{t}^{(j)}: t\ge 0)$ generated by the founder particle at the birth times $t_j$, $j=1,2,\ldots,N(t)$ (here we understand the sum is zero if $N(t)=0$). Let $X_{t}(x)$ be the number of particles alive at time $t$ with remaining lifetimes less than $x$. It is easy to see that $X_{t}(\infty):=\lim _{x \rightarrow \infty} X_{t}(x)=X_{t}((0, \infty))$. Let $H(s,t)=\mathbf{P}[s^{X_{t}(\infty)}]$ for $s\in [0,1]$ and $t\ge 0$. Then by (\ref{4.1'}) and (\ref{2.14}) we have
	\begin{align}\label{2.17}
		\!H(s,t)\!
		&=s^{1_{\{t-L<0\}}} \prod_{j=1}^{N(t)} H(s, t-t_j)\nonumber\\
		&=\int_0^{\infty}\! \mathbf{E}\Big[\prod_{j=1}^{N(t)} H(s, t-t_j) \Big| L \in d x\Big] s^{1_{\{t<x\}}} G(d x)\nonumber\\
		&=\int_0^{\infty}\! \mathbf{E}\Big[\exp \Big\{\int_0^{x\wedge t} \log H(s, t-r) N_x(dr)\Big\}\Big] s^{1_{\{t<x\}}} G(d x)\nonumber\\
		&=\int_0^{\infty}\!\! \exp \Big\{\!-\!\!\int_0^{x \wedge t} \!\sum_{n \in \mathbb{N}} \alpha(x-r) p(x-r, n)\big(1-e^{n \log H(s, t-r)}\big) d r\Big\} s^{1_{\{t<x\}}} G(d x)\nonumber\\
		&=\int_0^{\infty}\!\! \exp \Big\{\int_0^{x\wedge t} \!\alpha(x-r)\big[g(x-r, H(s, t-r))-1\big] d r\Big\} s^{1_{\{t<x\}}} G(d x).
	\end{align}
	It is obvious that the generation sizes of $(X_{t}: t\ge 0)$ form a Galton-Watson process. Let $\hat{\xi}_{n}$ be the number of particles in the $n$-th generation, i.e., $\hat{\xi}_{0}:=1$, $\hat{\xi}_{1}$ is the total number of children ever born to the ancestor, $\hat{\xi}_{2}$ is the total number of children ever born to members of the first generation, etc. Then we have
	\begin{proposition}\label{prop2.6}
		$(\widehat{\xi}_n: n \in \mathbb{N})$ is a Galton-Watson process with generating function
		$$
		h(s)=\int_0^{\infty} \exp \Big\{(s-1) \int_0^x \alpha(x-r) g'(x-r,1-) d r\Big\} G(d x),\quad s\in [0,1].
		$$
	\end{proposition}
	
	\begin{prof}
		Clearly the total number of offspring of different particles are independent and identically distributed random variables. For any $x\in (0,\infty)$,\ let $\lambda_{x}(\infty)=\mathbf{E}[N_{x}(\infty)]=\int_{0}^{x}\alpha(x-r)g'(x-r,1-)dr$. Then we have
		\begin{align*}
			\mathbf{P}(N_{x}(\infty)=n)=e^{-\lambda_{x}(\infty)}\frac{\lambda_{x}(\infty)^{n}}{n!},\quad x\in(0,\infty), n\in\mathbb{N}.
		\end{align*}
		So $\hat{\xi}_1=N(\infty)$ implies that
		\begin{align*}
			h(s)=\mathbf{E}\big[s^{\hat{\xi}_1}\big]
			&=\int_0^{\infty} \sum_{k=0}^{\infty} s^k \mathbf{P}[N_x(\infty)=k] G(d x)\\
			&=\int_0^{\infty}\Big[\sum_{k=0}^{\infty} e^{-s \lambda_x(\infty)} \frac{(s \lambda_x(\infty))^k}{k !}\Big] e^{(s-1) \lambda_x(\infty)} G(d x)\\
			&=\int_0^{\infty} \exp \Big\{(s-1) \int_0^x \alpha(x-r) g'(x-r,1-) d r\Big\} G(d x).
		\end{align*}
		Therefore it follows that
		\begin{align*}
			\mathbf{E}\big[s^{\hat{\xi}_2}\big]
			&=\sum_{k=0}^{\infty} \mathbf{E}\big[s^{\hat{\xi}_2} \big| \hat{\xi}_1=k\big] \mathbf{P}\big[\hat{\xi}_1=k\big]
			=\sum_{k=0}^{\infty} \mathbf{E}\big[s^{\sum_{i=1}^k N^{(i)}(\infty)}\big] \mathbf{P}[N(\infty)=k]\\
			&=\sum_{k=0}^{\infty} h(s)^k \mathbf{P}[N(\infty)=k]
			=h(h(s)),
		\end{align*}
		where the third equality holds since $\{N^{(i)}(\infty), i=1,2,\ldots,k\}$ and $N(\infty)$ are independent and identically distributed. In the same way we can check that $\mathbf{E}\big[s^{\hat{\xi}_n}\big]=h\big(\mathbf{E}\big[s^{\hat{\xi}_{n-1}}\big]\big)$, which is the characteristic of the Galton-Watson process.
		$\hfill\square$ 
	\end{prof}
	
	As a useful application of (\ref{2.17}), we now show that the extinction probability of $(X_{t}: t\ge 0)$ is the same as it of $(\widehat{\xi}_n: n \in \mathbb{N})$.
	
	\begin{theorem}\label{th2.7}
		Let $q(t)=\mathbf{P}[X_t (\infty)=0]$ and $q=\mathbf{P}[\lim_{n\rightarrow\infty}\hat{\xi}_{n}=0]$. Then $q(t) \uparrow q \in(0,1]$, as $t \rightarrow \infty$.
		In particular, $q=1$ if and only if  $\mathbf{E}[N(\infty)]\le 1$.
	\end{theorem}
	\begin{prof}
		Since $q(t)=H(0,t)$, it follows from (\ref{2.17}) that
		\begin{align}\label{2.18}
			q(t)=\int_0^t \exp \Big\{\int_0^x \alpha(x-r)[g(x-r, q(t-r))-1] d r\Big\} G(d x).
		\end{align}
		Notice that for any $t'\ge t$, $X_t (\infty)=0$ implies $X_{t'}(\infty)=0$. Then $q(t)\uparrow q_0$ as $t\rightarrow\infty$. Now it suffices to show $q_0 =q$. 
		Letting $t\rightarrow\infty$ in (\ref{2.18}) we get
		\begin{align*}
			q_0&=\int_0^{\infty} \exp \Big\{\int_0^x \alpha(x-r)[g(x-r, q_0)-1] d r\Big\} G(d x)\\
			&\ge\int_0^{\infty} \exp \Big\{\int_0^x \alpha(x-r)(q_0 -1)g'(x-r, 1-) d r\Big\} G(d x)=h(q_0)
		\end{align*}
		by monotone convergence. By the definition of $q$ we have $q$ is the smallest root in $[0,1]$ of $h(s)=s$. Thus the convexity of $h(s)$ implies $q_0 \ge q$. On the other hand, suppose that $\lim_{n\rightarrow\infty}\hat{\xi}_{n}\neq 0$, then infinitely many particles are realized. Therefore $\lim_{t\rightarrow\infty}{X_{t}(\infty)^{(+)}}=\infty$, where $X_{t}(\infty)^{(+)}$ is the number of particles which are born in $[0,t]$ (including the ancestor). If $\lim _{t \rightarrow \infty}X_t(\infty)=0$,  then $t_0 =\inf\{t>0: X_t(\infty)=0\}<\infty$ and $X_t(\infty)=0$ for $t\ge t_0$. Therefore for $t\ge t_0$ we have $X_{t}(\infty)^{(+)}=X_{t_0}(\infty)^{(+)}<\infty$. Hence
		\begin{align*}
			\Big\{\lim_{t \rightarrow \infty}X_t (\infty)=0\Big\}\subset\Big\{\lim_{t \rightarrow \infty}\hat{\xi}_n =0\Big\},
		\end{align*}
		which means that $q_0 \le q$. In particular, 
		\begin{align*}
			h'(s)\big|_{s=1}=\int_{0}^{\infty} G(dx) \int_{0}^{x} \alpha(x-s) g'(x-s,1-) ds=\mathbf{E}[N(\infty)].
		\end{align*}
		By the properties of the Galton-Watson process we have $q=1$ if and only if $\mathbf{E}[N(\infty)]\le 1$.
		$\hfill\square$ 
	\end{prof}

	\section{Existence of the non-degenerate limit}\label{section3}
	
	In this section we consider the sufficient and necessary condition for  $W_t^f:=e^{-\tilde{\alpha}t}\langle X_t,f\rangle$, $f\in B(0,\infty)^{+}$ to have a non-degenerate limit distribution as $t\rightarrow\infty$. Furthermore, the above convergence can be strengthened to hold with probability one under stronger sufficient conditions. We will discuss this case in Section \ref{section5}.
	
	If $W_t^f$ converges to some $W_{\infty}^f$ in distribution, then it follows from (\ref{2.3}), (\ref{2.5}) and (\ref{2.5'}) that $\phi^f(\theta)=\mathbf{E}\big[e^{-\theta W_\infty^f}\big]$ satisfies
	\begin{align}\label{3.6}
		\phi^f(\theta)=\int_0^{\infty} \exp \Big\{\int_0^x \alpha(x-s)\big[g(x-s, \phi^f(\theta e^{-\tilde{\alpha} s}))-1\big] d s\Big\} G(d x),\quad f\in B(0,\infty)^{+}.
	\end{align}
	Indeed, by (\ref{2.3}), as $t\rightarrow\infty$ we have
	\begin{align*}
		\mathbf{E}\big[e^{-\theta W_t^f}\big]=\mathbf{E}\big[e^{-\theta e^{-\tilde{\alpha} t}\langle X_t, f\rangle}\big]=\mathbf{E}\big[e^{-\langle X_t, \theta e^{-\tilde{\alpha}t} f\rangle}\big]=\big\langle G, e^{-u_t(\theta e^{-\tilde{\alpha}t} f)}\big\rangle \rightarrow \phi^f(\theta).
	\end{align*}
	By (\ref{2.5}) we have
		\begin{align}\label{3.1}
			&e^{-u_t(\theta e^{-\tilde{\alpha} t} f)(x)}-e^{-\theta e^{-\tilde{\alpha} t}f(x-t)}\nonumber\\
			&\quad=\int_0^t \alpha(x-s)\big[g\big(x-s,\langle G, e^{-u_{t-s}(\theta e^{-\tilde{\alpha}t} f)}\rangle\big)-1\big] e^{-u_{t-s}(\theta e^{-\tilde{\alpha}t} f)(x-s)} d s.
		\end{align}
	Notice that as $t\rightarrow\infty$, 
	\begin{align*}
		\langle G, e^{-u_{t-s}(\theta e^{-\tilde{\alpha}t} f)}\rangle=\langle G, e^{-u_{t-s}(\theta e^{-\tilde{\alpha}s}\cdot e^{-\tilde{\alpha}(t-s)} f)}\rangle=\mathbf{E}\big[e^{-\theta e^{-\tilde{\alpha}s} W_{t-s}^f}\big]\rightarrow \phi^f(\theta e^{-\tilde{\alpha}s}),
	\end{align*}
	and
	\begin{align*}
		u_t(\theta e^{-\tilde{\alpha}t} f)(x)
		&=\theta e^{-\tilde{\alpha}t} f(x-t)+\int_0^t \alpha(x-s)\big[1-g(x-s,\langle G, e^{-u_{t-s}(\theta e^{-\tilde{\alpha}t} f)}\rangle)\big] d s\\
		&\rightarrow \int_0^x \alpha(x-s)\big[1-g(x-s,\phi^f(\theta e^{-\tilde{\alpha}s}))\big] d s.
	\end{align*}
	Then letting $t\rightarrow\infty$ in (\ref{3.1}) we obtain (\ref{3.6}). Furthermore, if $\mathbf{E}W_{\infty}^f=a(f)$ given as (\ref{2.*'}), then $\phi^f (\theta)$ will also satisfy
	\begin{align}\label{3.7}
		\text{(A).}\ &0\le \phi^f (\theta)\le 1\ \text{for}\ \theta\in (0,\infty),\ \text{and}\ \phi^f (0)=1;\nonumber\\
		\text{(B).}\ &\phi^f\ \text{is continuous on}\ [0,\infty);\\
		\text{(C).}\ &\frac{1-\phi^f(\theta)}{\theta}\rightarrow a(f)\ \text{as}\ \theta\downarrow 0.\nonumber 
	\end{align}
	Then it is easy to investigate the uniqueness of solutions to  (\ref{3.6}) satisfying (\ref{3.7}).
	\begin{proposition}\label{prop3.1}
		Fix $f\in B(0,\infty)^{+}$. Suppose that $\phi_1^f$ and $\phi_2^f$ are solutions of (\ref{3.6}) which satisfy (\ref{3.7}). Then $\phi_1^f \equiv \phi_2^f$.
	\end{proposition}
	\begin{prof}
		Let $\tilde{\phi}^{f}(\theta)=\theta^{-1}\big|\phi_1^f(\theta)-\phi_2^f(\theta)\big|$ for $\theta>0$. Then it follows from (\ref{3.7}C) that $\tilde{\phi}^{f}(0+)=0$.\ Notice that  $|e^{-x_1}-e^{-x_2}| \le|x_1-x_2|$ for $x_1,x_2\ge 0$. Using the inequality we see from (\ref{3.6}) that for $\theta>0$,
		\begin{align*}
			\tilde{\phi}^{f}(\theta)
			&\le \theta^{-1} \int_0^{\infty}\Big|\exp \Big\{\int_0^x \alpha(x-s)\big[g\big(x-s, \phi_1^f(\theta e^{-\tilde{\alpha} s})\big)-1\big] d s\Big\}\\
			&\qquad\qquad\qquad-\exp \Big\{\int_0^x \alpha(x-s)\big[g\big(x-s, \phi_2^f(\theta e^{-\tilde{\alpha} s})\big)-1\big] d s\Big\}\Big| G(d x)\\
			&\le \theta^{-1} \int_0^{\infty}\Big|\int_0^x\!\! \alpha(x-s)\Big[g\big(x-s, \phi_2^f(\theta e^{-\tilde{\alpha} s})\big)-g\big(x-s, \phi_1^f(\theta e^{-\tilde{\alpha} s})\big)\Big] d s\Big| G(d x)\\
			&\le \int_0^{\infty}G(d x) \int_0^x \alpha(x-s) g'(x-s,1-)\cdot \theta^{-1} \Big|\phi_2^f(\theta e^{-\tilde{\alpha} s})-\phi_1^f(\theta e^{-\tilde{\alpha} s})\Big| d s\\
			&=\int_0^{\infty} \tilde{\phi}^f(\theta e^{-\tilde{\alpha} s})e^{-\tilde\alpha s} F(ds)\\
			&=\mathbf{E}\big[\tilde{\phi}^f(\theta e^{-\tilde{\alpha} \tilde{X}})\big],
		\end{align*}
		where $\tilde{X}$ is a random variable with $\mathbf{P}[\tilde{X}\le t]=\int_0^{t} e^{-\tilde\alpha s} F(ds)$. Iterating the above inequality we get $\tilde{\phi}^f(\theta) \leq \lim_{n \rightarrow \infty} \mathbf{E}\big[\tilde{\phi}^f(\theta e^{-\tilde\alpha S_n})\big]$, where $S_n$ is the sum of $n$ independent copies of $\tilde{X}$. Since $\mathbf{E}[\tilde{X}]>0$, by the strong large number law we have $e^{-\tilde\alpha S_n} \xrightarrow{\mathbf{P}\text{-a.s.}} 0$ as $n\rightarrow\infty$. Then by dominated convergence theorem we have $\tilde{\phi}^{f}(\theta)\le \tilde{\phi}^{f}(0+)=0$, which implies $\phi_1^f(\theta)=\phi_2^f(\theta)$ for $\theta>0$.
		$\hfill\square$ 
	\end{prof}
	
	We first define a random variable $Y$ which plays a very important role in studying the asymptotic behavior of $W_t^f$. Let
	\begin{align}\label{3.2}
		Y=\int_{0}^{\infty} e^{-\tilde{\alpha}s} N(ds),
	\end{align}
	where $N(t)$ is given as in Section \ref{section2}. Then it is simple to check the following proposition.
	
	\begin{proposition}\label{prop3.0}
		$\mathbf{E}Y=1$.
	\end{proposition}
	\begin{prof}
		Notice that $Y=\lim_{t\rightarrow\infty}\int_{0}^{t} e^{-\tilde{\alpha}s} N(ds)$. Then
		\begin{align*}
			\mathbf{E} Y=\lim _{t \rightarrow \infty}\Big\{\int_0^t \mathbf{E}\Big[\int_0^x e^{-\tilde{\alpha}s} N_x(d s)\Big] G(d x)+\int_t^{\infty} \mathbf{E}\Big[\int_0^t e^{-\tilde{\alpha} s} N_x(d s)\Big] G(d x)\Big\}.
		\end{align*}
		Since for any $x\in (0,\infty)$, the process $s \mapsto N_{x}(s)$ has at most countably many jumps, by It$\mathrm{\hat{o}}$'s formula, we have
		\begin{align*}
			\mathbf{E}\Big[\int_0^x e^{-\tilde{\alpha} s} N_x(d s)\Big]
			&=\mathbf{E}\Big[\int_0^x \int_0^{\alpha(x-s)} \int_{\mathbb{N}} \int_0^{p(x-s, n)} e^{-\tilde{\alpha} s} n M(d s, d u, d n, d v)\Big]\\
			&=\int_0^x e^{-\tilde{\alpha} s} \alpha(x-s) \sum_{n \in \mathbb{N}} p(x-s, n) n\ d s\\
			&=\int_0^x e^{-\tilde{\alpha} s} \alpha(x-s) g'(x-s, 1-) d s,\quad x\in (0,t].
		\end{align*}
		By similar calculations we also get
		\begin{align*}
			\mathbf{E}\Big[\int_0^t e^{-\tilde{\alpha} s} N_x(d s)\Big]=\int_0^t e^{-\tilde{\alpha} s} \alpha(x-s) g'(x-s, 1-) d s,\quad x\in [t,\infty).
		\end{align*}
		Then it is easy to check that
		\begin{align*}
			\mathbf{E}Y
			&=\lim_{t \rightarrow \infty}\int_{0}^{\infty} G(dx) \int_0^{x\wedge t} e^{-\tilde{\alpha} s} \alpha(x-s) g'(x-s, 1-) d s\\
			&=\lim_{t \rightarrow \infty}\int_0^t e^{-\tilde{\alpha} s} d s \int_s^{\infty} \alpha(x-s) g'(x-s, 1-) G(d x)=\int_0^{\infty} e^{-\tilde{\alpha} s} F(d s)=1,
		\end{align*}
		where the last equality follows by (\ref{2.*}).
		$\hfill\square$ 
	\end{prof}
	
	The key to our sufficient and necessary condition for the existence of a non-degenerate limit of $W_t^f$ or of a solution to (\ref{3.6}) satisfying (\ref{3.7}) lies in the following proposition. 
	
	\begin{proposition}\label{prop3.2}
		Let $Y$ be as in (\ref{3.2}) and
		\begin{align}\label{3.8}
			\psi(u)=u^{-1} \mathbf{E}\Big[e^{-u Y}+\int_0^{\infty}(1-\exp\{-u e^{-\tilde\alpha x}\}) N(d x)-1\Big],\quad u>0.
		\end{align}
		Then for any $0<r<1$, there exists some $\delta>0$ such that $\sum_{n=0}^{\infty} \psi(\delta r^n)<\infty$ if and only if
		\begin{align}\label{3.3}
			\mathbf{E}[Y|\log Y|]<\infty.
		\end{align}
	\end{proposition}
	
	\begin{prof}
		Since $\mathbf{E}Y=1$ and $e^{-x}$ is convex on $(0, \infty)$, by (\ref{3.8}) it is easy to see that $\psi(u)\ge 0$ for $u>0$. Notice that
		\begin{align*}
			\psi(u)
			&=u^{-1}\Big\{\int_0^{\infty} e^{\tilde{\alpha} x}[1-\exp \{-u e^{-\tilde{\alpha} x}\}] e^{-\tilde{\alpha} x} F(d x)-1+\mathbf{E}[e^{-uY}\}]\Big\}\\
			&=u^{-1}\{\mathbf{E}[Z^{-1}(1-e^{-u Z})]-1+\mathbf{E}[e^{-u Y}]\},\ 
		\end{align*}
		where $Z=e^{-\tilde{\alpha}\tilde{X}}$ with $\tilde{X}$ given as in the proof of Proposition \ref{prop3.1}. It follows that
		\begin{align*}
			\frac{d \psi(u)}{d u}
			&=u^{-2} \mathbf{E}\{Z^{-1}(u Z e^{-u Z}+e^{-u Z}-1)-(u Y e^{-u Y}+e^{-u Y}-1)\}\\
			&\rightarrow \frac{1}{2}[\mathbf{E}Y^2-\mathbf{E}Z]>0,\quad \text{as}\ u\downarrow 0.
		\end{align*}
		Then the continuity of $u\mapsto \frac{d \psi(u)}{d u}$ implies that $u\mapsto\psi(u)$ increases on $[0,u_1]$ for some $0<u_1<\infty$. By the integral test for convergence of positive series, for any $0<r<1$ and $0<\delta<\infty$ we have
		\begin{align*}
			\sum_{n=0}^{\infty} \psi(\delta r^n)<\infty \quad\text{if and only if}\quad \int_{0}^{\infty} \psi(\delta r^u)du<\infty.
		\end{align*}
		Notice that $\int_{0}^{\infty} \psi(\delta r^u)du=-(\log r)^{-1}\int_0^{\delta} u^{-1}\psi(u) du$. If $\int_0^{\delta} u^{-1}\psi(u) du<\infty$ holds for some $\delta>0$, it holds for all $\delta>0$. Therefore for any $0<r<1$ we have
		\begin{align*}
			\sum_{n=0}^{\infty} \psi(\delta r^n)<\infty\ \text{for some $\delta>0$} \quad\text{if and only if}\quad \int_0^{\delta} u^{-1}\psi(u) du<\infty\ \text{for all $\delta>0$}.
		\end{align*}
		On the other hand, it follows from Athreya \cite[Lemma 1]{Athreya69} that
		\begin{align*}
			\mathbf{E}[Y|\log Y|]<\infty\quad \text{if and only if}\quad \int_0^\delta u^{-2} \mathbf{E}[e^{-u Y}-e^{-u}] d u<\infty\ \text{for all $\delta>0$}.
		\end{align*}
		Let $\tilde{\psi}(u)=u^{-1} \mathbf{E}[e^{-uY}-e^{-u}] \ge 0$ for $u>0$. Then it suffices to show that for any $\delta>0$,
		\begin{align*}
			\int_0^{\delta} u^{-1}\psi(u) du<\infty\quad \text{if and only if}\quad \int_0^{\delta} u^{-1}\tilde\psi(u) du<\infty.
		\end{align*}
		To this end, we only need to show that for any $\delta>0$,
		\begin{align*}
			\int_0^{\delta} u^{-1}|\tilde\psi(u)-\psi(u)| du<\infty.
		\end{align*}
		Indeed, since $0 \le e^{-u}-1+u \leq \frac{1}{2} u^2$, it is easy to check that
		\begin{align*}
			0\le |\tilde\psi(u)-\psi(u)|
			&\le u^{-1}\{\mathbf{E}[Z^{-1}(e^{-u Z}-1+u Z)+(e^{-u}-1+u)]\}\\
			&\le u^{-1}\Big\{\frac{1}{2} \mathbf{E}[u^2 Z]+\frac{1}{2} u^2\Big\}\le u.
		\end{align*}
		Then the proposition follows. 
		$\hfill\square$ 
	\end{prof}
	
	We now show that (\ref{3.3}) is a necessary condition for the existence of a solution to (\ref{3.6}) satisfying (\ref{3.7}A), (\ref{3.7}B) and a more general version of  (\ref{3.7}C).
	
	\begin{proposition}\label{prop3.3}
		Fix $f\in B(0,\infty)^{+}$. Let $0<c(f)<\infty$. Suppose that there exists a unique solution to (\ref{3.6}) satisfying (\ref{3.7}A), (\ref{3.7}B) and
		\begin{align}\label{3.7'}
			\frac{1-\phi^f(\theta)}{\theta}\rightarrow c(f),\quad \text{as}\ \theta\downarrow 0.
		\end{align} 
		Then (\ref{3.3}) holds.
	\end{proposition}
	
	\begin{prof}
		Suppose that there exists a unique solution to (\ref{3.6}) satisfying (\ref{3.7}A), (\ref{3.7}B) and (\ref{3.7'}). Let
		\begin{align*}
			b_1^f(\theta)=\frac{1-\phi^f(\theta)}{\theta},\quad \theta>0.
		\end{align*}
		Then $b_1^f(\theta)\ge 0$ and  $b_1^f(0+)=c(f)$. It follows that there exists some positive constants $c_1, c_2, c_3<\infty$, for any $\theta\in(0,c_1]$, $c_2\le b_1^f(\theta)\le c_3$. Then for any $\theta\in(0,c_1]$ we have $0 \le \phi^f(\theta) \le 1-c_2 \theta \le e^{-c_2 \theta} \le 1$. Notice that
			\begin{align*}
				\mathbf{E}[e^{-\theta Y}]
				&=\lim _{t \rightarrow \infty}\int_0^{\infty}\! \mathbf{E}\Big[\exp\Big\{\!-\!\!\int_0^{x\wedge t} \theta e^{-\tilde{\alpha}s} N_x(d s)\Big\}\Big] G(d x)\\
				&=\lim _{t \rightarrow \infty}\int_0^{\infty}\! \mathbf{E}\Big[\exp \Big\{\!-\!\!\int_0^{x\wedge t}\!\!\! \int_0^{\alpha(x-s)}\!\!\! \int_{\mathbb{N}} \int_0^{p(x-s, n)}\!\!\! n \theta e^{-\tilde{\alpha} s} M(d s, d u, d n, d v)\Big\}\Big] G(d x)\\
				&=\int_0^{\infty} \exp \Big\{-\int_0^x \alpha(x-s) \sum_{n \in \mathbb{N}} p(x-s, n)\big(1-e^{-n\theta e^{-\tilde\alpha s}}\big) d s\Big\} G(d x)\\
				&=\int_0^{\infty} \exp \Big\{\int_0^x \alpha(x-s)\big[g(x-s, e^{-\theta e^{-\tilde{\alpha} s}})-1\big] d s\Big\} G(d x).
			\end{align*}
		For any $\theta\in(0,c_1]$, by (\ref{3.6}) we have
			\begin{align}\label{3.9}
				b_1^f(\theta)
				=&\ \theta^{-1}\Big[1-\int_0^{\infty} \exp \Big\{\int_0^x \alpha(x-s)\big[g(x-s, \phi^f(\theta e^{-\tilde{\alpha} s}))-1\big] d s\Big\} G(d x)\Big]\nonumber\\
				=&\ -\theta^{-1}\Big[\int_0^{\infty} \exp \Big\{\int_0^x \alpha(x-s)\big[g(x-s, \phi^f(\theta e^{-\tilde{\alpha} s}))-1\big] d s\Big\} G(d x)\nonumber\\
				&\qquad\qquad+\int_0^{\infty}\big[1-\phi^f(\theta e^{-\tilde{\alpha} s})\big] F(d s)-1\Big]+\theta^{-1} \int_0^{\infty}\big[1-\phi^f(\theta e^{-\tilde{\alpha} s})\big] F(d s)\nonumber\\
				\le&\ -\theta^{-1}\Big[\int_0^{\infty} \exp \Big\{\int_0^x \alpha(x-s)\big[g(x-s, e^{-c_2\theta e^{-\tilde{\alpha} s}})-1\big] d s\Big\} G(d x)\nonumber\\
				&\qquad\qquad+\int_0^{\infty}\big[1-e^{-c_2\theta e^{-\tilde{\alpha} s}}\big] F(d s)-1\Big]+\int_0^{\infty}b_1^f(\theta e^{-\tilde{\alpha}s}) e^{-\tilde{\alpha}s} F(d s)\nonumber\\				
				=&\ \int_0^{\infty}b_1^f(\theta e^{-\tilde{\alpha}s}) e^{-\tilde{\alpha}s} F(d s)-c_2 \psi(c_2 \theta),
			\end{align}
		where $\psi$ is given by (\ref{3.8}). Since $t\mapsto F(t)$ is continuous, there exists $T_0>0$ with $1<c_4:=\int_0^{T_0} e^{-\tilde{\alpha}s} F(ds)<1$. Let $b_2^f (\theta)=\sup_{0\le x\le \theta} b_1^f(x)$. Then it follows from (\ref{3.9}) that for any $\theta\in (0,c_1]$ we have
		\begin{align*}
			b_2^f (\theta)
			&\le \sup_{0\le x\le \theta}\!\Big\{\int_0^{\infty}\!b_1^f(x e^{-\tilde{\alpha}s}) e^{-\tilde{\alpha}s} F(d s)\Big\}-c_2 \psi(c_2 \theta)\\
			&= \sup_{0\le x\le \theta}\!\Big\{\int_0^{T_0}\!b_1^f(x e^{-\tilde{\alpha}s}) e^{-\tilde{\alpha}s} F(d s)\Big\}+\!\sup_{0\le x\le \theta}\!\Big\{\!\int_{T_0}^{\infty}\!b_1^f(x e^{-\tilde{\alpha}s}) e^{-\tilde{\alpha}s} F(d s)\Big\}-c_2 \psi(c_2 \theta).
		\end{align*}
		For any $s\in (0,T_0]$, we have $e^{-\tilde{\alpha}s}\in (0,1]$, which implies that $0\le xe^{-\tilde{\alpha}s}\le \theta \le c_1$. Then it follows from $b_1^f(\cdot)\le 0$ that  
		\begin{align*}
			\sup_{0\le x\le \theta}\Big\{\int_0^{T_0}b_1^f(x e^{-\tilde{\alpha}s}) e^{-\tilde{\alpha}s} F(d s)\Big\}\le b_2^f(\theta) \int_0^{T_0} e^{-\tilde\alpha s} F(d s)=c_4 b_2^f(\theta).
		\end{align*}
		On the other hand, since $\theta\mapsto b_2^f(\cdot)\le 0$ increases,\ it is natural to see that
		\begin{align*}
			\sup_{0\le x\le \theta}\Big\{\int_{T_0}^{\infty}b_1^f(x e^{-\tilde{\alpha}s}) e^{-\tilde{\alpha}s} F(d s)\Big\}
			&\le \int_{T_0}^{\infty}b_2^f(\theta e^{-\tilde{\alpha}s}) e^{-\tilde{\alpha}s} F(d s)\\
			&\le b_2^f(\theta e^{-\tilde{\alpha}T_0})\Big[\int_0^{\infty} e^{-\tilde\alpha s} F(d s)-c_4\Big]\\
			&=(1-c_4)b_2^f(\theta e^{-\tilde{\alpha}T_0}).
		\end{align*}
		Hence we have
		\begin{align*}
			b_2^f (\theta)
			&\le c_4 b_2^f(\theta)+(1-c_4)b_2^f(\theta e^{-\tilde{\alpha}T_0})-c_2 \psi(c_2 \theta)\\
			&\le b_2^f(r\theta)-c_5 \psi(c_2 \theta),
		\end{align*}
		where $r=e^{-\tilde\alpha T_0} \in(0,1),\ c_5=c_2 /\left(1-c_4\right)$. Iterating the above inequality yields
		\begin{align*}
			b_2^f(\theta) \le b_2^f(r^{m+1} \theta)-c_5 \sum_{n=0}^m \psi(c_2 \theta r^n).
		\end{align*}
		Notice that $b_2^f(\theta)$ is bounded on $(0,c_1]$, then the above inequality implies $\sum_{n=0}^{\infty} \psi(c_2 \theta r^n)<\infty$. Then by Proposition \ref{prop3.2} we get (\ref{3.3}). 
		$\hfill\square$ 
	\end{prof}
	
	We next show that (\ref{3.3}) implies that  $\mathbf{E}[e^{-\theta W_t^f}]\rightarrow \phi^f(\theta)$ as $t\rightarrow\infty$, where $\phi^f(\theta)$ is the unique solution of (\ref{3.6}) satisfying (\ref{3.7}). We need the following preliminary proposition.
	
	\begin{proposition}\label{prop3.4}
		Fix $f\in B(0,\infty)^{+}$. Let $I^f(\theta, t)=\theta^{-1} \mathbf{E}\big[e^{-\theta W_t^f}+\theta W_t^f-1\big]$. Suppose that (\ref{3.3}) holds, then
		\begin{align*}
			\lim _{\theta \downarrow 0} \sup _{t\ge 0}|I^f(\theta, t)|=0.
		\end{align*}
	\end{proposition}
	
	\begin{prof}
		Let $b_3(x)=x^{-1}(e^{-x}+x-1)$ for $x>0$. Then $b_3(x)\ge 0$, $x\mapsto b_3(x)$ increases and $b_3(0+)=0$. It follows that $I^f(\theta, t)=\mathbf{E}\big[W_t^f b_3(\theta W_t^f)\big] \ge 0$. Let $I_T^f(\theta)=\sup _{0\le t \le T} I^f(\theta, t)$ for $T \in(0, \infty)$. Thus for any fixed $T$ we have $I_T^f(\theta)=\mathbf{E}\big[W_T^f b_3(\theta W_T^f)\big]$, where $W_T^f=\sup _{0 \le t \le T} W_t^f<\infty$, $\mathbf{P}$-a.s. by Cheng and Li \cite[Proposition 3.4]{Cheng23+}. This means that $\theta\mapsto I_T^f(\theta)$ increases and $I_T^f(0+)=0$ for any fixed $T$. It follows from $I^f(\theta, t)\ge 0$ that $\mathbf{E}\big[e^{-\theta W_t^f}\big]\ge 1-a(f)\theta$, where $a(f)$ is given by (\ref{2.*'}). Then there exists some positive constants $c_6, c_7<\infty$, for any $t\ge 0$ and $\theta\in(0,c_7]$ we have $\mathbf{E}\big[e^{-\theta W_t^f}\big]\ge 1-a(f)\theta\ge e^{-c_6 \theta}$. Notice that
		\begin{align}\label{3.**}
			\mathbf{E}\big[e^{-\theta W_t^f}\big]
			&=\int_0^{\infty} \exp \big\{-u_t(\theta e^{-\tilde\alpha t} f)(x)\big\} G(d x)\nonumber\\
			&=\int_0^{\infty}\!e^{-\theta e^{-\tilde{\alpha} t}f(x-t)}\cdot\exp\Big\{\int_0^t  \alpha(x-s)\big[g\big(x-s,\langle G, e^{-u_{t-s}(\theta e^{-\tilde{\alpha}t} f)}\rangle\big)-1\big]  d s\Big\}G(dx)\nonumber\\
			&=\!\int_t^{\infty}\!\!\! \big(e^{-\theta e^{-\tilde\alpha t} f(x-t)}\!-\!1\big)\exp \Big\{\!\int_0^t\! \alpha(x-s)\big[g\big(x-s, \mathbf{E}\big[e^{-\theta e^{-\tilde\alpha s}W_{t-s}^f}\big]\big)-1\big] d s\Big\} G(d x)\nonumber\\			
			&\quad+\int_0^{\infty} \exp \Big\{\int_0^t \alpha(x-s)\big[g\big(x-s, \mathbf{E}\big[e^{-\theta e^{-\tilde\alpha s}W_{t-s}^f}\big]\big)-1\big] d s\Big\} G(d x)
		\end{align}
		and
		\begin{align*}
			\mathbf{E} W_t^f=\int_0^t \mathbf{E} W_{t-s}^f\ e^{-\tilde{\alpha} s} F(d s)+e^{-\tilde{\alpha} t} \int_t^{\infty} f(x-t) G(d x).
		\end{align*}
		For $u>0$ and $t\ge 0$, let
		\begin{align}\label{3.10}
			\Psi(u, t)=u^{-1} \mathbf{E}\Big[e^{-u \int_0^t e^{-\tilde\alpha s} N(d s)}+\int_0^t\big(1-e^{-u e^{-\tilde\alpha s}}\big) N(d s)-1\Big].
		\end{align}
		It is easy to see that
		\begin{align*}
			\Psi(u, t)=u^{-1}\Big\{
			&\int_0^{t} \exp \Big\{\int_0^x \alpha(x-s)\big[g\big(x-s, e^{-u e^{-\tilde\alpha s}}\big)-1\big] d s\Big\} G(d x)\\
			&+\int_t^{\infty} \exp \Big\{\int_0^t \alpha(x-s)\big[g\big(x-s, e^{-u e^{-\tilde\alpha s}}\big)-1\big] d s\Big\} G(d x)\\
			&+\int_0^t\big(1-e^{-u e^{-\tilde\alpha s}}\big) d s \int_s^{\infty} \alpha(x-s) g'(x-s, 1-) G(d x)-1\Big\}.
		\end{align*}
		Then for any $u>0$ we have $\frac{\partial \Psi(u, t)}{\partial t}\ge 0$. Thus it follows from $\lim _{t \rightarrow 0+} \Psi(u, t)=0$ that for any $u>0$ we have $\Psi(u, t)\ge 0$ and $t\mapsto \Psi(u, t)$ increases. Moreover, for any $u>0$ we have $\lim _{t \rightarrow \infty} \Psi(u, t)=\psi(u)$, where $\psi$ is given by (\ref{3.8}). Then for any $\theta\in (0,c_7]$,
		\begin{align*}
			&I^f(\theta, t)\\
			=&\ \frac{1}{\theta}\Big\{\mathbf{E}\big[e^{-\theta W_t^f}\big]+\theta \int_0^t \mathbf{E} W_{t-s}^f\ e^{-\tilde{\alpha}s} F(d s)+\theta e^{-\tilde\alpha t} \int_t^{\infty} f(x-t) G(d x)-1\Big\}\\
			=&\ \int_0^t\! \mathbf{E} W_{t-s}^f e^{-\tilde{\alpha} s} F(d s)-1+\!\!\int_0^{\infty}\!\!\! \exp \Big\{\int_0^t\! \alpha(x-s)\big[g(x-s, \mathbf{E}[e^{-\theta e^{-\tilde\alpha s} W_{t-s}^f}])-1\big] d s\Big\} G(d x)\\
			&\ +\frac{1}{\theta}\!\int_t^{\infty}\!\!\!\!\big(1\!-\!e^{-\theta e^{-\tilde\alpha t}f(x-t)}\big)\Big[1\!-\!\exp \Big\{\!\int_0^t\!\! \alpha(x-s)\big[g(x-s, \mathbf{E}[e^{-\theta e^{-\tilde\alpha s} W_{t-s}^f}])-1\big] d s\Big\}\Big]G(d x)\\
			&\ +\frac{1}{\theta}\int_t^{\infty}\Big[e^{-\theta e^{-\tilde{\alpha} t} f(x-t)}+\theta e^{-\tilde\alpha t} f(x-t)-1\Big] G(d x)\\
			\le&\ \frac{1}{\theta}\Big\{\int_0^t \mathbf{E}\big[e^{-\theta e^{-\tilde\alpha s} W_{t-s}^f}\big] e^{-\tilde\alpha s} F(d s)+\theta \int_0^t\! \mathbf{E} W_{t-s}^f e^{-\tilde{\alpha} s} F(d s)-\int_0^t e^{-\tilde\alpha s}F(ds)\Big\}\\
			&\ +\frac{1}{\theta}\Big\{\int_0^{\infty}\exp \Big\{\int_0^t \alpha(x-s)\big[g(x-s, \mathbf{E}[e^{-\theta e^{-\tilde\alpha s} W_{t-s}^f}])-1\big] d s\Big\}G(d x)-1\\
			&\qquad\quad+\int_0^t \mathbf{E}\big[1-e^{-\theta e^{-\tilde\alpha s} W_{t-s}^f}\big] e^{-\tilde\alpha s} F(d s)\Big\}\\
			&\ +\int_t^{\infty} e^{-\tilde{\alpha} t} f(x-t) G(d x)\int_0^t \alpha(x-s)\big[1-g(x-s, \mathbf{E}[e^{-\theta e^{-\tilde\alpha s} W_{t-s}^f}])\big] d s\\
			&\ +\theta \int_t^{\infty} \frac{1}{2} e^{-2\tilde\alpha t}(f(x-t))^2 G(d x)\\
			\le&\ \int_0^t I^f(\theta e^{-\tilde\alpha s}, t-s) e^{-\tilde\alpha s} F(d s)\\
			&\ +\!\frac{1}{\theta}\Big\{\!\int_0^{\infty}\!\!\!\!\!\exp \Big\{\!\int_0^t\!\! \alpha(x-s)\big[g(x-s, e^{-c_6 \theta e^{-\tilde\alpha s}})\!-\!1\big] d s\Big\}G(d x)\!-\!1\!+\!\!\int_0^t\!\! \big[1\!-\!e^{-c_6 \theta e^{-\tilde\alpha s}}\big] F(d s)\Big\}\\
			&\ +\|f\| \int_0^{\infty} G(d x) \int_0^t \alpha(x-s) g'(x-s, 1-) c_6 \theta e^{-\tilde\alpha s} d s+\frac{1}{2}\|f\| \theta\\
			=&\ \int_0^t I^f(\theta e^{-\tilde\alpha s}, t-s) e^{-\tilde\alpha s} F(d s)+c_6 \Psi(c_6 \theta,t) +\|f\|\int_0^t c_6 \theta e^{-\tilde\alpha s}F(ds)+\frac{1}{2}\|f\| \theta\\
			\le&\ \int_0^t I^f(\theta e^{-\tilde\alpha s}, t-s) e^{-\tilde\alpha s} F(d s)+c_6 \psi(c_6 \theta) +c_8 \theta,
		\end{align*}
		where $c_8 =\|f\|(\frac{1}{2}+c_6)\in (0,\infty)$. This means that
		\begin{align}\label{3.11}
			I_T(\theta)\le \int_0^t I_T^f(\theta e^{-\tilde\alpha s}) e^{-\tilde\alpha s} F(d s)+c_6 \psi(c_6 \theta) +c_8 \theta.
		\end{align}
		By arguments similar to those in the proof of Proposition \ref{prop3.3} we have
		\begin{align*}
			\int_0^T I_T^f\left(\theta e^{-\tilde\alpha s}\right) e^{-\tilde\alpha s} F(d s) \le c_4 I_T^f(\theta)+(1-c_4) I_T^f(\theta e^{-\tilde\alpha T_0}),
		\end{align*}
		where $c_4$ and $T_0$ are as in the proof of Proposition \ref{prop3.3}. Then by (\ref{3.11}) we have
		\begin{align*}
			I_T^f(\theta) \le I_T^f(\theta e^{-\tilde\alpha T_0})+(1-c_4)^{-1}[c_6 \psi(c_6 \theta)+c_8 \theta].
		\end{align*}
		Since $I_T^f(0+)=0$, iterating the above inequality yields
		\begin{align}\label{3.12}
			I_T^f(\theta) \le(1-c_4)^{-1}\Big[c_6 \sum_{n=0}^{\infty} \psi(c_6 \theta r^n)+c_8 \theta /(1-e^{-\tilde\alpha T_0})\Big],
		\end{align}
		where $r=e^{-\tilde{\alpha} T_0}\in (0,1)$. Notice that the right-hand side of the above inequality is independent of $T$ and is finite by Proposition \ref{prop3.2}. Since $u\mapsto\psi(u)$ increases and $\psi(0+)=0$, by letting $T\rightarrow\infty$ first and then $\theta\downarrow 0$ we get the desired result.
		$\hfill\square$ 
	\end{prof}
	
	\begin{proposition}\label{prop3.5}
		Fix $f\in B(0,\infty)^{+}$. Suppose that (\ref{3.3}) holds. Then $\lim_{t\rightarrow\infty}\mathbf{E}[e^{-\theta W_t^f}]= \phi^f(\theta)$, where $\phi^f(\theta)$ is the unique solution of (\ref{3.6}) satisfying (\ref{3.7}).
	\end{proposition}
	
	\begin{prof}
		For any $\theta>0$ and $T\ge 0$, let $J^f(\theta, t)=\theta^{-1}\big\{\mathbf{E}\big[e^{-\theta W_t^f}\big]-\phi^f(\theta)\big\}$, $J_T^f(\theta)=\sup _{t\ge T}|J^f(\theta, t)|$, $J^f(\theta)=\lim _{T \rightarrow \infty} J_T^f(\theta)$. Then it suffices to show that for any $\theta\in (0,\infty)$ we have $J^f(\theta) \equiv 0$. By the definitions of $I^f (\theta,t)$ and $J^f(\theta,t)$ we have 
		\begin{align*}
			|J^f(\theta,t)-I^f(\theta,t)|
			&=\big|\theta^{-1}\big\{\theta \mathbf{E}W_t^f +\phi^f(\theta)-1\big\}\big|\\
			&\le \big|\mathbf{E} W_t^f-a(f)\big|+\big|\theta^{-1}(1-\phi^f(\theta))-a(f)\big|.
		\end{align*}
		It follows from (\ref{3.7}C) and Proposition \ref{prop3.4} that $J^f (0+)=0$. Let
			\begin{align*}
				b(\theta,f,t)&=\int_0^t \Big[\exp\Big\{\int_0^x \alpha(x-s)\big[g(x-s, \mathbf{E}[e^{-\theta e^{-\tilde\alpha s} W_{2t-s}^f}])-1\big] d s\Big\}\\
				&\qquad\quad-\exp\Big\{\int_0^x \alpha(x-s)\big[g(x-s, \phi^f(\theta e^{-\tilde\alpha s}))-1\big] d s\Big\}\Big]G(dx).
			\end{align*}
		Using (\ref{3.**}) and (\ref{3.6}) we obtain
			\begin{align}\label{3.13}
				&\big|\theta J^f(\theta,2t)-b(\theta,f,t)\big|\nonumber\\
				=&\ \Big|\int_t^{2t}\exp\Big\{\int_0^x \alpha(x-s)\big[g(x-s, \mathbf{E}[e^{-\theta e^{-\tilde\alpha s} W_{2t-s}^f}])-1\big] d s\Big\} G(dx)\nonumber\\
				&\ +\!\!\int_{2t}^{\infty}\!\!e^{-\theta e^{-2\tilde\alpha t}f(x-2t)}\exp\Big\{\!\int_0^{2t}\!\!\! \alpha(x-s)\big[g(x-s, \mathbf{E}[e^{-\theta e^{-\tilde\alpha s} W_{2t-s}^f}])-1\big] d s\Big\} G(dx)\Big|\nonumber\\
				\le&\int_t^{2 t} G(d x)+\int_{2 t}^{\infty} G(d x)=1-G(t).
			\end{align}
		Notice that
			\begin{align*}
				|b(\theta,f,t)|
				&\le\int_0^{t}G(dx)\int_0^x \alpha(x-s)\Big|g(x-s, \mathbf{E}[e^{-\theta e^{-\tilde\alpha s} W_{2t-s}^f}])-g(x-s, \phi^f(\theta e^{-\tilde\alpha s}))\Big| d s\\
				&\le\int_0^{t}G(dx)\int_0^t \alpha(x-s)g'(x-s,1-)\Big|\mathbf{E}[e^{-\theta e^{-\tilde\alpha s} W_{2t-s}^f}]-\phi^f(\theta e^{-\tilde\alpha s})\Big| d s\\
				&=\theta \int_0^{\infty}\big|J^f(\theta e^{-\tilde{\alpha} s} , 2 t-s)\big| e^{-\tilde{\alpha} s} F(d s).
			\end{align*}
		Then by (\ref{3.13}) we have
		\begin{align*}
			J_{2 T}^f(\theta)
			&=\sup _{t \ge T} J^f(\theta, 2 t)\\
			&\le \sup _{t \ge T} \big\{\theta^{-1}|b(\theta,f,t)|+\theta^{-1}[1-G(t)]\big\}\\
			&\le \int_0^{\infty}\big|J_T^f(\theta e^{-\tilde{\alpha} s}, 2 t-s)\big| e^{-\tilde{\alpha} s} F(d s)+\theta^{-1}[1-G(T)]\\
			&=\mathbf{E}\big[J_T^f(\theta e^{-\tilde{\alpha} \tilde{X}})\big]+\theta^{-1}[1-G(T)],
		\end{align*}
		where $\tilde{X}$ is given as in the proof of Proposition \ref{prop3.1}. By letting $T\rightarrow\infty$ we get $J^f(\theta)\le\mathbf{E}\big[J^f(\theta e^{-\tilde{\alpha} \tilde{X}})\big]$. Iterating the above inequality we get $J^f(\theta)\le\lim_{n \rightarrow \infty}\mathbf{E}\big[J^f(\theta e^{-\tilde{\alpha} S_n})\big]$, where $S_n$ is the sum of $n$ independent copies of $\tilde{X}$. Since $\mathbf{E}[\tilde{X}]>0$, by the strong large number law we have $e^{-\tilde\alpha S_n} \xrightarrow{\mathbf{P}\text{-a.s.}} 0$ as $n\rightarrow\infty$. Then by dominated convergence theorem we have $J^{f}(\theta)\le J^{f}(0+)=0$. Hence it follows from the non-negativity of $J^f(\theta)$ that $J^f(\theta) \equiv 0$, for $\theta \in(0, \infty)$.
		$\hfill\square$ 
	\end{prof}
	
	\begin{theorem}\label{th3.A}
		Suppose that (\ref{2.*}) and (\ref{super}) hold. Then for any $f\in B(0,\infty)^+$, as $t\rightarrow \infty$, $W_t^f \xrightarrow{d} W_{\infty}^{f}$ exists and is not identically zero if and only if (\ref{3.3}) holds. In this case, we have
		\begin{align*}
			\rm{(1).}&\ \mathbf{E} W_{\infty}^f=a(f)\  \text{given by (\ref{2.*'})};\\
			\rm{(2).}&\ \mathbf{P}[W_\infty ^{f}=0]=q\  \text{(the extinction probability)};\\
			\rm{(3).}&\ \phi^f(\theta)=\mathbf{E}\big[e^{-\theta W_\infty^f}\big]\ \text{is the unique solution of (\ref{3.6}) satisfying (\ref{3.7}).}
		\end{align*}
	\end{theorem}
	
	\begin{prof}
		Suppose that (\ref{3.3}) holds. For any fixed $f\in B(0,\infty)^+$, then by Proposition \ref{prop3.5} we have $\lim_{t\rightarrow\infty}\mathbf{E}[e^{-\theta W_t^f}]= \phi^f(\theta)$, where $\phi^f(\theta)$ is the unique solution of (\ref{3.6}) satisfying (\ref{3.7}). It follows that $\phi^f(\theta)$ is continuous and $\phi^f(\theta)=\mathbf{E}[e^{-\theta W_\infty^f}]$. Thus by (\ref{3.7}C) we have $\mathbf{E} W_{\infty}^f=a(f)$. Recall that
		\begin{align*}
			q=\mathbf{P}\big[\lim _{t \rightarrow \infty}X_t(\infty)=0\big]=\int_0^{\infty} \exp \Big\{\int_0^x \alpha(x-r)[g(x-r, q)-1] d r\Big\} G(d x)<1.
		\end{align*}
		Let $q^*:=\mathbf{P}\big[W_{\infty}^f=0\big]=\lim _{\theta \rightarrow \infty} \phi^f(\theta)$. Then by (\ref{3.6}) we have
		\begin{align*}
			q^*=\int_0^{\infty} \exp \Big\{\int_0^x \alpha(x-r)[g(x-r, q^*)-1] d r\Big\} G(d x)\in [0,1].
		\end{align*}
		Notice that
		\begin{align}\label{3.***}
			b_4(y):=\int_0^{\infty} \exp \Big\{\int_0^x \alpha(x-r)[g(x-r, y)-1] d r\Big\} G(d x)
		\end{align}
		is a convex function on $[0,1]$ with  $b'_4(0)<1$ and $b'_4(1)>1$. The only two fixed points of $b_4(y)$ on $[0,1]$ are $q$ and $1$, hence that $q^*=q$ or $q^*=1$. Since $q^*=1$ implies $\mathbf{E} W_{\infty}^f=0$, we must have $q^*=q$. Conversely, suppose that $W_{\infty}^{f}=\lim_{t\rightarrow \infty} W_t^f$ exists and is not identically zero for any $f\in B(0,\infty)^+$. Then $\lim _{t \rightarrow \infty} \mathbf{E} W_t^f=a(f)$ implies $0<\mathbf{E} W_{\infty}^f=c (f)\le a(f)$ by Fatou's lemma. Since by (\ref{2.4}) and (\ref{2.5'}) we have $\phi^f(\theta)=\mathbf{E}\big[e^{-\theta W_\infty^f}\big]$ is the unique solution of (\ref{3.6}) satisfying (\ref{3.7}A), (\ref{3.7}B) and (\ref{3.7'}). Then it follows from Proposition \ref{prop3.1} that (\ref{3.3}) holds.
		$\hfill\square$ 
	\end{prof}

	\section{Convergence of the Malthusian normalized random measures}\label{section4}
	
	By Theorem \ref{th3.A} we have $W_t^{1}=e^{-\tilde{\alpha}t} X_t(\infty)\xrightarrow{d}W_\infty^1$ exists and is not identically zero if (\ref{2.*}), (\ref{super}) and (\ref{3.3}) hold. In this section we want to express $W_{\infty}^f$ as the product of a nonrandom functional $f\mapsto A(f)$ and  $W_\infty^1$ for any $f\in B(0,\infty)^+$. Then we could get the convergence of the Malthusian normalized random measures $e^{-\tilde{\alpha}t} X_t$. To this end, We just need to discuss the convergence in probability with respect to $\mathbf{P}$ of the age distribution $A_t(f):=\langle X_t,f\rangle/X_t(\infty)$ as $t\rightarrow\infty$. To ensure $A_t(f)$ is well defined a.e., we establish the results conditioning on the non-extinction event.
	
	We just need to consider the case where $f(y)=1_{(0,x]}(y)$, $x,y\in (0,\infty)$. 
	Let $A_t(x):=X_t(x)/X_t(\infty)$ for $x\in (0,\infty)$. For the convenience of statement of the results, in the rest of the paper, we write $(X_t^y:t\ge 0)$ be the process defined on $(\Omega, \mathcal{F}, \mathcal{F}_{t}, \mathbf{P})$ with the same distribution as $\hat{X}=(\Omega, \mathcal{F}, \mathcal{F}_{t}, \hat{X}_{t}, \hat{\mathbf{P}}_{y})$ for any fixed $y\in(0,\infty)$, i.e., 
	\begin{align}\label{4.y}
		\mathbf{P}\big[X_t^{y}\in\cdot\big]=\hat{\mathbf{P}}_{y}\big[\hat{X}_t\in\cdot\big],\quad y\in (0,\infty).
	\end{align} 
	We start with a simple but useful equality about $X_t(x)$, which is followed from the additive property of branching processes. For any $x\in(0,\infty)$ we write
	\begin{align}\label{4.2}
		X_{t+s}(x)=\sum_{i=1}^{X_t(\infty)} X_s^{x_i}(x),\quad t,s\ge 0,
	\end{align}
	where $\{x_i; i=1,2,\ldots,X_t(\infty)\}$ is the remaining-lifetime chart at time $t$, $\big\{X_s^{x_i}; i=1,2,\ldots,X_t(\infty)\big\}$ are independent and further if $x_i=y$ then the conditional distribution of $X_s^{x_i}$ is the same as $X_s^{y}$ defined as (\ref{4.y}).
	Then $\mathbf{E}[X_s^{x_i}(x)]=\pi 1_{(0,x]}(x_i)$. It follows from the above equality (\ref{4.2}) that
	\begin{align}\label{4.3}
		e^{-\tilde{\alpha} s} \frac{X_{t+s}(x)}{X_t(\infty)}
		=\ &\frac{1}{X_t(\infty)} \sum_{i=1}^{X_t(\infty)}\big[X_s^{x_i}(x)-\pi_s 1_{(0, x]}(x_i)\big] e^{-\tilde{\alpha} s}\nonumber\\
		&+\frac{1}{X_t(\infty)} \sum_{i=1}^{X_t(\infty)}\big[\pi_s 1_{(0, x]}(x_i) e^{-\tilde{\alpha} s}-n_1 V(x_i)A(x)\big]\nonumber\\
		&+\frac{\langle X_t,V\rangle}{X_t(\infty)} n_1 A(x)\nonumber\\
		=:\ &a_t(x,s)+b_t(x,s)+c_t A(x),
	\end{align}
	where
	\begin{align}\label{4.n}
		n_1=\frac{\int_0^{\infty} e^{-\tilde{\alpha} u} d u \int_u^{\infty} 1_{(0, \infty)}(y-u) G(d y)}{\int_0^{\infty} u e^{-\tilde{\alpha} u} F(d u)}=\frac{\int_0^{\infty} e^{-\tilde{\alpha} u}[1-G(u)] d u}{\int_0^{\infty} u e^{-\tilde{\alpha} u} F(d u)},
	\end{align}
	\begin{align}\label{4.a}
		A(x)=\frac{\int_0^{\infty} e^{-\tilde\alpha u} d u \int_u^{\infty} 1_{(0, x]}(y-u) G(d y)}{\int_0^{\infty} e^{-\tilde\alpha u} d u \int_u^{\infty} 1_{(0, \infty)}(y-u) G(d y)}
	\end{align}
	and
	\begin{align}\label{4.v}
		V(x)=\int_0^{\infty} \alpha(x-s) g'(x-s, 1-) e^{-\tilde\alpha s} d s.
	\end{align}
	Then it is easy to check that
	\begin{align}\label{4.4}
		A_{t+s}(x)=\frac{a_t(x, s)+b_t(x, s)+c_t A(x)}{a_t(\infty, s)+b_t(\infty, s)+c_t}.
	\end{align}
	
	We first show that for any $t\ge 0$ and $x\in (0,\infty)$, as $s\rightarrow\infty$ we have $|b_t(x,s)|\xrightarrow{\mathbf{P}\text{-a.s.}}0$ and $|b_t(\infty,s)|\xrightarrow{\mathbf{P}\text{-a.s.}}0$.
	
	\begin{proposition}\label{prop4.1}
		Let $n_1$, $A(x)$ and $V(x)$ be given by (\ref{4.n}), (\ref{4.a}) and (\ref{4.v}), respectively. Then for any $y\in (0,\infty)$, as $s\rightarrow\infty$ we have
		\begin{align*}
			\big|\pi_s 1_{(0, x]}(y) e^{-\tilde{\alpha} s}-n_1 V(y)A(x)\big|\rightarrow 0\quad\text{and}\quad \big|\pi_s 1_{(0, \infty)}(y) e^{-\tilde{\alpha} s}-n_1 V(y)\big|\rightarrow 0.
		\end{align*}
	\end{proposition}
	
	\begin{prof}
		It follows from (\ref{2.2'}) that
		\begin{align*}
			\pi_{s} 1_{(0,x]}(y)=1_{(0,x]}(y-s)+\int_{0}^{s} \alpha(y-r) g'(y-r,1-) \langle G,\pi_{s-r}1_{(0,x]}\rangle dr.
		\end{align*}
		By integrating both sides of the above equality with respect to $G(dx)$ we obtain
		\begin{align*}
			\langle G,\pi_{s} 1_{(0,x]}\rangle
			&=G(x+s)-G(s)+\int_{0}^{s} \langle G,\pi_{s-r}1_{(0,x]}\rangle F(dr).
		\end{align*}
		Notice that for any $x\in (0,\infty)$, by Proposition \ref{prop2.2} we have
		\begin{align*}
			\lim_{t \rightarrow \infty}\langle G, \pi_t 1_{(0, x]}\rangle e^{-\tilde{\alpha} t} =a(1_{(0, x]})=n_1 A(x).
		\end{align*}
		Then for any $x\in (0,\infty)$ and $\varepsilon>0$, there exists $T>0$ such that for $t\ge T$, we have $e^{-\tilde\alpha t}<\varepsilon$ and $\big|\langle G, \pi_t 1_{(0, x]}\rangle e^{-\tilde{\alpha} t}-n_1 A(x)\big|<\varepsilon$. Hence for $s>[1+1\vee \beta\alpha^{-1}]T$, by Cheng and Li \cite[Proposition 3.4]{Cheng23+} we have
		\begin{align*}
			&\big|\pi_s 1_{(0, x]}(y) e^{-\tilde{\alpha} s}-n_1 V(y)A(x)\big|\\
			=&\ \Big|1_{(0, x]}(y-s) e^{\tilde\alpha s}+e^{-\tilde\alpha s}\! \int_0^s\! \alpha(y-r) g'(y-r, 1-)\langle G, \pi_{s-r} 1_{(0, x]}\rangle d r-n_1 A(x) V(y)\Big|\\
			\le&\ e^{-\tilde{\alpha} s}+\int_0^{s-T} \alpha(y-r) g'(y-r, 1-)\big|\langle G, \pi_{s-r} 1_{(0, x]}\rangle e^{-\tilde{\alpha}(s-r)}-n_1 A(x)\big| e^{-\tilde{\alpha} r} d r\\
			\quad&+\int_{s-T}^s \alpha(y-r) g^{\prime}\left(y-r,1-\right)\langle G, \pi_{s-r} 1_{(0, x]}\rangle e^{-\tilde\alpha s} d r\\
			\quad&+n_1 A(x) \int_{s-T}^{\infty} \alpha(y-r) g'(y-r, 1-) e^{-\tilde{\alpha} r} d r\\
			<&\ [2+\beta\alpha^{-1}(1+n_1 A(x))]\varepsilon,
		\end{align*}
		which follows that $\big|\pi_s 1_{(0, x]}(y) e^{-\tilde{\alpha} s}-n_1 V(y)A(x)\big|\rightarrow 0$ as $s\rightarrow\infty$. On the other hand, notice that by Proposition \ref{prop2.2} we also have
		\begin{align*}
			\lim_{t \rightarrow \infty}\langle G, \pi_t 1_{(0, \infty)}\rangle e^{-\tilde{\alpha} t} =a(1_{(0, x]})=n_1.
		\end{align*}
		Then the second result can be obtained in a similar way.
		$\hfill\square$ 
	\end{prof}
	
	Next we prove the following two results conditioning on the non-extinction event:
	\begin{description}
		\item[\rm{(i).}] For fixed $s$, $a_t(x,s)\xrightarrow{\mathbf{P}} 0$ and $a_t(\infty,s)\xrightarrow{\mathbf{P}} 0$ as $t\rightarrow\infty$;
		
		\item[\rm{(ii).}] $c_t$ is bounded below in probability.
	\end{description}
	
	\begin{proposition}\label{prop4.2}
		For any fixed $s\in (0,\infty)$, conditioning on the non-extinction event, we have
		\begin{align*}
			\frac{1}{X_t(\infty)} \sum_{i=1}^{X_t(\infty)}\big[X_s^{x_i}(x)-\pi_s 1_{(0, x]}(x_i)\big]\xrightarrow{\mathbf{P}} 0,\quad \text{as}\ t\rightarrow\infty.
		\end{align*}
	\end{proposition}
	
	\begin{prof}
		Notice that $X_s^{x_i}(x)$ are nonnegative random variables and $\sup_{x,y}\pi_s 1_{(0,x]}(y)\le e^{\beta s}<\infty$ by Cheng and Li \cite[Proposition 3.4]{Cheng23+}. Then it suffices to show that for any $0<\theta<\infty$ we have
		\begin{align*}
			\mathbf{E}\Big[\exp\Big\{-\frac{\theta}{X_t(\infty)} \sum_{i=1}^{X_t(\infty)}\big[X_s^{x_i}(x)-\pi_s 1_{(0, x]}(x_i)\big] \Big\}\Big]\rightarrow 1,\quad\text{as}\ t\rightarrow\infty.
		\end{align*}
		It is simple to see that
		\begin{align}\label{4.4'}
			\mathbf{E}&\Big[\exp\Big\{-\frac{\theta}{X_t(\infty)} \sum_{i=1}^{X_t(\infty)}\big[X_s^{x_i}(x)-\pi_s 1_{(0, x]}(x_i)\big] \Big\}\Big|\mathcal{F}_t\Big]\nonumber\\
			&=\exp\Big\{\sum_{i=1}^{X_t(\infty)}\Big[\frac{\theta}{X_t(\infty)}\pi_s 1_{(0, x]}(x_i)+\log \mathbf{E}\big(e^{-\frac{\theta}{X_t(\infty)}X_s^{x_i}(x)}\big)\Big]\Big\}.
		\end{align}
		Using the facts $\log(1-x)=-x+o(x)$ as $x\rightarrow 0$ and 
		\begin{align*}
			\Big|\frac{1-\mathbf{E}[e^{-\theta_1 X_s^{y}(x)}]}{\theta_1}-\pi_s 1_{(0, x]}(y)\Big| \rightarrow 0,\quad \text{as}\ \theta_1\downarrow 0,
		\end{align*}
		for any $x,y\in (0,\infty)$ we have
		\begin{align*}
			\Big|\pi_s 1_{(0, x]}(y)+\frac{X_t(\infty)}{\theta} \log \mathbf{E}\big[e^{-\frac{\theta}{X_t(\infty)} X_s^y(x)}\big]\Big| \rightarrow 0,\quad \text {as}\ X_t(\infty) \rightarrow \infty .
		\end{align*}
		Then by (\ref{4.4'}) we obtain the desired result.
		$\hfill\square$ 
	\end{prof}
	
	\begin{proposition}\label{prop4.3}
		Suppose that (\ref{2.*}), (\ref{super}) and (\ref{3.3}) hold. Then for any $\varepsilon>0$, conditioning on the non-extinction event, there exists $\delta>0$ such that
		\begin{align*}
			\liminf_{t\rightarrow\infty}\mathbf{P}\Big[\frac{\langle X_t,V\rangle}{X_t(\infty)}>\delta\Big]>1-\varepsilon.
		\end{align*} 
	\end{proposition}
	
	\begin{prof}
		It follows from Theorem \ref{th3.A} that 
		\begin{align*}
			e^{-\tilde{\alpha}t}X_t(\infty)\xrightarrow{d}W_\infty^1,\ e^{-\tilde{\alpha}t}\langle X_t,V\rangle\xrightarrow{d}W_\infty^V,\quad\text{as}\ t\rightarrow\infty,
		\end{align*}
		and (\ref{3.3}) implies that conditioning on the non-extinction event we have $\mathbf{P}[W_\infty^1>0]=\mathbf{P}[W_\infty^V>0]=1$. Then the proposition now follows easily.
		$\hfill\square$ 
	\end{prof}
	
	Thus we now establish the following results as the consequences of Propositions \ref{prop4.1}-\ref{prop4.3} and Theorem \ref{th3.A}.
	
	\begin{theorem}\label{th4.2}
		Suppose that (\ref{2.*}), (\ref{super}) and (\ref{3.3}) hold. Then for any $f\in B(0,\infty)^+$, conditioning on the non-extinction event, as $t\rightarrow\infty$ we have
		\begin{align}\label{4.*a}
			A_t(f)\xrightarrow{\mathbf{P}}A(f):=\frac{\int_0^{\infty} e^{-\tilde\alpha u} d u \int_u^{\infty} f(y-u) G(d y)}{\int_0^{\infty} e^{-\tilde\alpha u} [1-G(u)]d u} <\infty.
		\end{align}
	\end{theorem}
	
	\begin{prof}
		Using (\ref{4.4}), Propositions \ref{prop4.1}-\ref{prop4.3} and the monotone convergence theorem we easily obtain the desired result.
		$\hfill\square$ 
	\end{prof}
	
	\begin{theorem}\label{th4.A}
		Suppose that (\ref{2.*}) and (\ref{super}) hold. Then for any $f\in B(0,\infty)^+$, as $t\rightarrow \infty$, 
		$$W_t^f \xrightarrow{d} A(f)W_{\infty}^{1}$$ 
		exists and is not identically zero if and only if (\ref{3.3}) holds, where $W_{\infty}^{1}$ given as in Theorem \ref{th3.A} is the limit (in distribution sense) of $W_{t}^{1}$. Furthermore, the distribution of $A(f)W_{\infty}^{1}$ is the same as $W_{\infty}^{f}$ given in Theorem \ref{th3.A}.
	\end{theorem}
	
	\begin{theorem}\label{th4.A'}
		Suppose that (\ref{2.*}) and (\ref{super}) hold. Then there is a finite measure $Q$ on $\mathfrak{M}(0,\infty)$ such that $Q(\mathbf{0})=q$ (the extinction probability) and
		\begin{align*}
			\mathbf{P}\big[e^{-\tilde{\alpha}t}X_t\in\cdot\big]\xrightarrow{w}Q(\cdot)
		\end{align*}
		if and only if (\ref{3.3}) holds, where $\mathbf{0}$ denotes the null measure and $\xrightarrow{w}$ stands for weak convergence. In this case, the Laplace transform of $Q$ is given by
		\begin{align*}
			\int_{\mathfrak{M}(0,\infty)} e^{-\langle \nu,f\rangle}Q(d\nu)=e^{-A(f)}\phi^1(1),
		\end{align*}
		where $\phi^1(1)$ is given as in Theorem \ref{th3.A}.
	\end{theorem}
	
	\begin{prof}
		It follows from Theorem \ref{th4.A} that as $t\rightarrow\infty$,
		\begin{align*}
			\int_{\mathfrak{M}(0,\infty)} e^{-\langle \nu,f\rangle}\mathbf{P}(e^{-\tilde{\alpha}t}X_t\in d\nu)=\mathbf{E}\big[e^{-W_t^f}\big]\rightarrow \mathbf{E}\big[e^{-A(f)W_{\infty}^1}\big]=e^{-A(f)}\phi^1(1).
		\end{align*}
		Notice that $A$ is a functional on $B(0,\ \infty)^{+}$continuous with respect to bounded pointwise convergence. Then the result follows by the continuity theorem for the Laplace function of random measures; see, e.g., Li \cite[Theorem 1.20]{Li22}.
		$\hfill\square$ 
	\end{prof}
	
	As an application of Theorem \ref{th4.A}, we consider the absolute continuity of $W_{\infty}^f$ given in Theorem \ref{th3.A}. We omit the proof of the following proposition since the argument is similar to Doney \cite[Theorem B]{Doney72b} or Doney \cite[Theorem  7.7]{Doney72a}.
	
	\begin{proposition}\label{thB}
		Suppose that (\ref{2.*}), (\ref{super}) and (\ref{3.3}) hold. Then there exists a continuous function $w(x)\ge 0$ such that
		\begin{align*}
			\mathbf{P}\big[x_1<W_{\infty}^{1}\le x_2\big]=\int_{x_1}^{x_2}w(x) dx,\quad\text{for}\   0<x_1<x_2\le\infty,
		\end{align*}
		where $W_{\infty}^{1}$ given as in Theorem \ref{th3.A} is the limit (in distribution sense) of $W_{t}^{1}$.
	\end{proposition}
	
	\begin{theorem}\label{thB'}
		Suppose that (\ref{2.*}), (\ref{super}) and (\ref{3.3}) hold. Then for any $f\in B(0,\infty)^+$,
		\begin{align*}
			\mathbf{P}\big[x_1<W_{\infty}^{f}\le x_2\big]=A(f)^{-1}\int_{x_1}^{x_2}w(A(f)^{-1} x) dx,\quad\text{for}\  0<x_1<x_2\le\infty,
		\end{align*}
		where $A(f)$ is given as (\ref{4.*a}) and $w(x)$ is given in Proposition \ref{thB}.
	\end{theorem}
	
	\begin{prof}
		Notice that the distribution of $W_{\infty}^{f}$ is the same as $A(f)W_{\infty}^{1}$ by Theorem \ref{th4.A}. Then for any $0<x_1<x_2\le\infty$, by Proposition \ref{thB} we have
		\begin{align*}
			\mathbf{P}\big[x_1<A(f) W_{\infty}^1 \le x_2\big]
			&=\mathbf{P}\big[A(f)^{-1}x_1<W_{\infty}^1 \le A(f)^{-1}x_2\big]\\
			&=\int_{A(f)^{-1}x_1}^{A(f)^{-1}x_2}w(x)dx=A(f)^{-1}\int_{x_1}^{x_2}w(A(f)^{-1} x) dx,
		\end{align*}
		which completes the proof.
		$\hfill\square$ 
	\end{prof}

	\section{Almost sure convergence}\label{section5}
	
	For any $f\in B(0,\infty)^+$, as $t\rightarrow\infty$, it follows from Theorem \ref{th4.A} that $W_t^{f}\xrightarrow{d}A(f)W_\infty^1$ exists and is not identically zero if (\ref{2.*}), (\ref{super}) and (\ref{3.3}) hold. In this section, we want to show that under some ``$L \log L$'' assumptions we have $W_t^f$ converges almost surely to a non-degenerate random variables with the same distribution as $A(f) W_{\infty}^1$. To this end, we need some preliminary results.
	
	Recall that $(X_t^y:t\ge 0)$ is the process defined on $(\Omega, \mathcal{F}, \mathcal{F}_{t}, \mathbf{P})$ with the same distribution as $\hat{X}=(\Omega, \mathcal{F}, \mathcal{F}_{t}, \hat{X}_{t}, \hat{\mathbf{P}}_{y})$ for any fixed $y\in(0,\infty)$. For any $y\in (0,\infty)$, we define $\tau_1(y)=\inf\{t>0:X_{t}^y\neq X_{t-}^y\}$ with the convention $\inf\emptyset=\infty$, which denotes the first branching time of $(X_t^y:t\ge 0)$. Let $\eta_{\tau_1(y)}$ be the number of offsprings produced by the ancestor at its remaining lifetime $y-\tau_1(y)$. Then by the definition of $(X_t^y:t\ge 0)$ we have
	\begin{align*}
		\mathbf{P}[\tau_1(y)\in ds]=\alpha(y-s)e^{-\int_{0}^{s}\alpha(y-r)dr}ds,\quad s>0,
	\end{align*}
	and
	\begin{align*}
		\mathbf{P}[\eta_{\tau_1(y)}=n]=p(y-\tau_1(y),n),\quad n\in \mathbb{N}.
	\end{align*}
	Then using the two random variables we study the $L\log L$-type moments of $(X_t:t\ge 0)$ and $(X_t^y:t\ge 0)$.
	
	\begin{theorem}\label{th5.*}
		Suppose that 
		\begin{align}\label{5.*}
			\sup _{y \ge 0} \sum_{n=0}^{\infty} n|\log n| p(y, n)<\infty.
		\end{align}
		Then for any $x\in (0,\infty]$ and $t\ge 0$ we have $\mathbf{E}[X_t(x)|\log X_t(x)|]<\infty$.
	\end{theorem}
	
	\begin{prof}
		Write $H_1(x):=x|\log x|$ and $\mu_t(x):=\mathbf{E}[H_1(X_t(x))]$. Observe that $H_1(x)$ is a nonnegative and convex function on $[1,\infty)$ and there exists $K>0$, for any $x,y\in [1,\infty)$ we have $H_1(x,y)\le K H_1(x)H_2(y)$. Notice that for any $x\in(0,\infty)$, $X_t(x)$ can be represented as the sum
		\begin{align}\label{5.0}
			X_t(x)=\sum_{j=1}^{\eta_{\tau_1(L)}} X_{t-\tau_1(L)}^{(j)}(x) + 1_{\{0<L-t\le x\}},
		\end{align}
		involving $\eta_{\tau_1(L)}$ independent daughter processes $(X_{t}^{(j)}: t\ge 0)$ generated by the founder particle at the birth times $\tau_1(L)$, $j=1,2,\ldots,\eta_{\tau_1(L)}$. Then we obtain that
			\begin{align*}
				\mu_t(x)\!
				&=\mathbf{E}[H_1(X_t(x))]\\
				&=\mathbf{E}\Big[H_1\Big(\sum_{j=1}^{\eta_{\tau_1(L)}} X_{t-\tau_1(L)}^{(j)}(x) + 1_{\{0<L-t\le x\}}\Big)\Big]\\
				&\le \frac{K}{2}H_1(2)\mathbf{E}\Big[H_1\Big(\sum_{j=1}^{\eta_{\tau_1(L)}} X_{t-\tau_1(L)}^{(j)}(x)\Big)\Big]+\frac{H_1(2)}{2}\\
				&\le \frac{K}{2}H_1(2)\int_{0}^{\infty}G(dy)\int_0^t \sum_{n=0}^{\infty} \mathbf{E}\Big[H_1\Big(\sum_{j=1}^{n} X_{t-s}^{(j)}(x)\Big)\Big]\mathbf{P}[\eta_{s}=n]\mathbf{P}[\tau_1(y)\in ds] +\frac{H_1(2)}{2}\\
				&\le \frac{K}{2}H_1(2)\!\!\int_{0}^{\infty}\!\!\!\!G(dy)\!\!\int_0^t \sum_{n=0}^{\infty} \Big\{\frac{1}{n}\sum_{j=1}^{n}\mathbf{E}\big[H_1\big(n X_{t-s}^{(j)}(x)\big)\big]\!\Big\}\mathbf{P}[\eta_{s}=n]\mathbf{P}[\tau_1(y)\in ds]\! +\!\frac{H_1(2)}{2}\\
				&\le \frac{K^2}{2}H_1(2)\int_{0}^{\infty}G(dy)\int_0^t \sum_{n=0}^{\infty} H_1(n)\mu_{t-s}(x)\mathbf{P}[\eta_{s}=n]\mathbf{P}[\tau_1(y)\in ds] +\frac{H_1(2)}{2}\\
				&\le \frac{K^2}{2}H_1(2)\!\int_{0}^{\infty}\!\!\!G(dy)\!\!\int_0^t \sum_{n=0}^{\infty} \!H_1(n)\mu_{t-s}(x)p(y-s, n) \alpha(y-s) e^{-\int_0^s \alpha(y-r) d r} d s +\frac{H_1(2)}{2}\\
				&\le \frac{K^2}{2}H_1(2)\|\alpha\|
				\cdot \sup _{y \ge 0} \sum_{n=0}^{\infty} H_1(n) p(y, n) \int_0^t \mu_{s}(x) d s +\frac{H_1(2)}{2}.
			\end{align*}
		Then the desired result follows from Gronwall's inequality.
		$\hfill\square$ 
	\end{prof}
	
	\begin{corollary}\label{cor5.*}
		Suppose that (\ref{5.*}) holds. 
		Then for any $x,y\in (0,\infty]$ and $t\ge 0$ we have $\mathbf{E}[X_t^y(x)|\log X_t^y(x)|]<\infty$.
	\end{corollary}
	
	\begin{prof}
		For any $x,y\in(0,\infty)$, it follows by arguments similar to (\ref{2.14}) that $X_t^y(x)$ can be represented as the sum
		\begin{align*}
			X_t^y(x)=\sum_{j=1}^{N_y(t)} X_{t-t_j}^{(j)}(x) + 1_{\{0<y-t\le x\}},
		\end{align*}
		where $N_y(t)$ is given as in (\ref{4.1'}) and $(X_{t}^{(j)}: t\ge 0)$,  $j=1,2,\ldots,N_y(t)$ are independent and have the same distribution as $(X_{t}: t\ge 0)$. Recall that $H_1(x)=x|\log x|$ and $\mu_t(x)=\mathbf{E}[H_1(X_t(x))]$. Then
		\begin{align*}
			\lambda_y(t):=\mathbf{E}[N_y(t)]=\int_0^{y \wedge t} \alpha(y-s) g'(y-s, 1-) d s<\infty.
		\end{align*}
		Fix $t\ge 0$, for any $0\le s\le t$ and $x\in (0,\infty)$, it follows from the proof of Theorem \ref{th5.*} that there exists $C(t)<\infty$ such that $\mu_s(x)\le C(t)$. Then we obtain that
		\begin{align*}
			\mathbf{E}[H_1(X_t^y(x))]
			&=\mathbf{E}\Big[H_1\Big(\sum_{j=1}^{N_y(t)} X_{t-t_j}^{(j)}(x) + 1_{\{0<y-t\le x\}}\Big)\Big]\\
			&\le \frac{K}{2}H_1(2)\mathbf{E}\Big[H_1\Big(\sum_{j=1}^{N_y(t)} X_{t-t_j}^{(j)}(x)\Big)\Big]+\frac{H_1(2)}{2}\\
			&\le \frac{K^2}{2}H_1(2)\mathbf{E}\Big[\frac{1}{N_y(t)}H_1(N_y(t))\sum_{j=1}^{N_y(t)}H_1 \big(X_{t-t_j}^{(j)}(x)\big)\Big]+\frac{H_1(2)}{2}\\
			&=\frac{K^2}{2}H_1(2)\sum_{n=1}^{\infty}\frac{H_1(n)}{n}\sum_{j=1}^{n}\mathbf{E}\Big[H_1 \big(X_{t-t_j}^{(j)}(x)\big)\Big]e^{\lambda_y(t)}\frac{(\lambda_y(t))^n}{n!}\!+\!\frac{H_1(2)}{2}\\
			&=\frac{K^2}{2}H_1(2)e^{\lambda_y(t)}\sum_{n=1}^{\infty}\log n \sum_{j=1}^{n}\mu_{t-t_j}(x)\frac{(\lambda_y(t))^n}{n!}+\frac{H_1(2)}{2}\\
			&\le \frac{K^2}{2}H_1(2)C(t)e^{\lambda_y(t)}\sum_{n=1}^{\infty}n \frac{(\lambda_y(t))^n}{(n-1)!}+\frac{H_1(2)}{2}\\
			&= \frac{K^2}{2}H_1(2)C(t)e^{\lambda_y(t)}\Big[\sum_{n=2}^{\infty} \frac{(\lambda_y(t))^n}{(n-2)!}+\sum_{n=1}^{\infty} \frac{(\lambda_y(t))^n}{(n-1)!}\Big]+\frac{H_1(2)}{2}\\
			&= \frac{K^2}{2}H_1(2)C(t)\lambda_y(t)\big[\lambda_y(t)+1\big]+\frac{H_1(2)}{2}\\
			&<\infty,
		\end{align*}
		which completes the proof.
		$\hfill\square$ 
	\end{prof}
	
	We now restrict the index $t$ to lattices of the form $\{n\delta; n\in \mathbb{N}, \delta\ \text{is a positive rational}\}$ and further consider the almost sure convergence of the age distribution $A_{n\delta} (f)=\langle X_{n\delta},f\rangle/X_{n\delta}(\infty)$ as $n\rightarrow\infty$ conditioning on the non-extinction event.
	
	\begin{proposition}\label{prop5.2}
		Let $V\in B(0,\infty)^+$ be given by (\ref{4.v}). Then the process $(\langle X_t ,V\rangle e^{-\tilde{\alpha}t}:t\ge 0)$ is a martingale with respect to the filtration $(\mathcal{F}_t)_{t\ge 0}$.
	\end{proposition}
	
	\begin{prof}
		It follows by arguments similar to (\ref{4.2}) that
		\begin{align*}
			\langle X_{t+s},V\rangle=\sum_{i=1}^{X_t(\infty)} \langle X_s^{x_i},V\rangle,\quad t,s\ge 0,
		\end{align*}
		where $\{x_i; i=1,2,\ldots,X_t(\infty)\}$ is the remaining-lifetime chart at time $t$, $\big\{X_s^{x_i}; i=1,2,\ldots,X_t(\infty)\big\}$ are independent and further if $x_i=y$ then the conditional distribution of $X_s^{x_i}$ is the same as $X_s^{y}$ defined as (\ref{4.y}). It suffices to show that for any $t\ge 0$ and $x\in (0,\infty)$ we have $\pi_t V(x)=V(x)e^{\tilde{\alpha}t}$. Indeed, by (\ref{2.2'}) we obtain
			\begin{align*}
				\pi_t V(x)
				&=V(x-t)+\int_0^t \alpha(x-s) g'(x-s, 1-)\langle G, \pi_{t-s} V\rangle d s\\
				&=\int_0^{\infty}\!\! \alpha(x-t-r) g'(x-t-r, 1-) e^{-\tilde{\alpha} r} d r+\!\int_0^t\! \alpha(x-s) g'(x-s, 1-)\langle G, \pi_{t-s} V\rangle d s\\
				&=V(x) e^{\tilde{\alpha}t}+\int_0^t\alpha(x-s) g'(x-s, 1-)\big[\langle G, \pi_{t-s} V\rangle-e^{\tilde{\alpha}(t-s)}\big] d s. 			
			\end{align*}
		Notice that
		\begin{align*}
			\langle G, V\rangle=\int_0^{\infty} G(d x) \int_0^{\infty} \alpha(x-r) g'(x-r, 1-) e^{-\tilde{\alpha} r} d r=\int_0^{\infty} e^{-\tilde{\alpha} r} F(d r)=1.
		\end{align*}
		Then
		\begin{align*}
			\big|\pi_t V(x)-V(x) e^{\tilde{\alpha}t}\big| \le \int_0^t \alpha(x-s) g^{\prime}(x-s, 1-)\big|\langle G, \pi_{t-s} V\rangle-\langle G, V\rangle e^{\tilde{\alpha}(t-s)} \big| d s,
		\end{align*}
		and thus
		\begin{align*}
			\big\|\pi_t V(\cdot)-V(\cdot) e^{\tilde{\alpha} t}\big\| \le \beta \int_0^t\big\|\pi_s V(\cdot)-V(\cdot) e^{\tilde\alpha s}\big\| d s.
		\end{align*}
		By Gronwall's inequality we have $\big\|\pi_t V(\cdot)-V(\cdot) e^{\tilde{\alpha} t}\big\|=0$, which means $\pi_t V(x)=V(x) e^{\tilde{\alpha} t}$ for any $x\in (0,\infty)$. Therefore, for any $s,  t\ge 0$ we obtain
		\begin{align*}
			\mathbf{E}\big[\langle X_{t+s}, V\rangle e^{-\tilde{\alpha}(t+s)} \big| \mathcal{F}_t\big]
			&=\mathbf{E}\Big[\sum_{i=1}^{X_t(\infty)}\langle X^{x_i}_{s}, V\rangle \Big| \mathcal{F}_t\Big]e^{-\tilde{\alpha}(t+s)}\\
			&=\sum_{i=1}^{X_t(\infty)}\hat{\mathbf{E}}_{x_i}\big[\langle \hat{X}_{s}, V\rangle \big]e^{-\tilde{\alpha}(t+s)}\\
			&=\sum_{i=1}^{X_t(\infty)}\pi_s V(x_i)e^{-\tilde{\alpha}(t+s)}\\
			&=\sum_{i=1}^{X_t(\infty)}
			V(x_i)e^{-\tilde{\alpha}t}\\
			&=\langle X_t,V \rangle e^{-\tilde{\alpha}t},
		\end{align*}
		which completes the proof.
		$\hfill\square$ 
	\end{prof}
	
	\begin{proposition}\label{prop5.2'}
		Suppose that (\ref{2.*}), (\ref{super}) and (\ref{5.*}) hold. Then for any $m\in\mathbb{N}$ and fixed $\delta>0$, conditioning on the non-extinction event, we have
		\begin{align}\label{5.5}
			\frac{1}{X_{n\delta}(\infty)} \sum_{i=1}^{X_{n\delta}(\infty)}\big[X_{m\delta}^{x_i}(x)-\pi_{m\delta} 1_{(0, x]}(x_i)\big]\xrightarrow{\mathbf{P}\text{-a.s.}} 0,\quad \text{as}\ n\rightarrow\infty.
		\end{align}		
	\end{proposition}
	
	\begin{prof}
		By Proposition \ref{prop2.2},  conditioning on the non-extinction event, there exists $0<\hat{\alpha}\le\tilde{\alpha}$ and constant $c_8>0$ such that $X_{n\delta}(\infty)\ge c_8 e^{\hat{\alpha}n\delta}$. It follows from Cheng and Li \cite[Proposition 3.4]{Cheng23+} that for any fixed $m\delta$, $\sup_y \pi_{m \delta} 1_{(0, x]}(y)<\infty$. Then for any $z>\sup_y \pi_{m \delta} 1_{(0, x]}(y)$ we have
		\begin{align*}
			\mathbf{P}\big[\big|X_{m \delta}^y(x)-\pi_{m\delta} 1_{(0 ,x]}(y)\big| \ge z\big] \le \mathbf{P}[X_{m\delta}^y(x) \ge z].
		\end{align*}
		By Corollary \ref{cor5.*} we have $\mathbf{E}\big[X_{m \delta}^y(x)\big|\log X_{m \delta}^y(x)\big|\big]<\infty$. Then for any $\varepsilon>0$, we conclude from Athreya \cite[Proposition 1]{Athreya76} that 
		\begin{align*}
			\sum_{n=1}^{\infty}\mathbf{P}\Big\{\Big|\frac{1}{X_{n\delta}(\infty)} \sum_{i=1}^{X_{n\delta}(\infty)}\big[X_{m\delta}^{x_i}(x)-\pi_{m\delta} 1_{(0, x]}(x_i)\big]\Big|>\varepsilon \ \Big|\mathcal{F}_{n\delta}\Big\}<\infty,\quad \mathbf{P}\text{-a.s.}.
		\end{align*}
		Hence we obtain the desired result by the extended Borel-Cantelli lemma; see, e.g., Breiman \cite[Proposition 5.29]{Breiman68}.
		$\hfill\square$ 
	\end{prof}
	
	\begin{proposition}\label{prop5.3}
		Suppose that (\ref{2.*}), (\ref{super}), (\ref{3.3}) and (\ref{5.*}) hold. Then for any fixed  $\delta>0$, conditioning on the non-extinction event, we have
		\begin{align*}
			\lim_{n\rightarrow\infty}\frac{\langle X_{n\delta},V\rangle}{X_{n\delta}(\infty)}>0,\quad \mathbf{P}\text{-a.s.}.
		\end{align*}
	\end{proposition}
	
	\begin{prof}
		It follows from Proposition \ref{prop5.2} that as $t\rightarrow\infty$, 
		\begin{align}\label{5.v}
			\langle X_{t},V\rangle e^{-\tilde{\alpha}t}\xrightarrow{\mathbf{P}\text{-a.s.}}\widetilde{W}_{\infty}^{V}
		\end{align} exists and further the distribution of $\widetilde{W}_{\infty}^{V}$ is the same as $W_{\infty}^{V}$ given as in Theorem \ref{th3.A}. Then conditioning on the non-extinction event, we have
		$\mathbf{P}[\widetilde{W}_{\infty}^{V}>0]=1$ and $C_V:=\sup_{t}\langle X_{t},V\rangle e^{-\tilde{\alpha}t}<\infty$ by Theorem \ref{th3.A}. Therefore it  suffices to show 
		\begin{align*}
			\limsup_{n \rightarrow \infty} e^{-\tilde{\alpha} n \delta} X_{n \delta}(\infty)<\infty,\quad  \mathbf{P}\text {-a.s.}.
		\end{align*}
		Let $0<\varepsilon<1/2$. There exists $n_0$ such that 
		\begin{align*}
			\sup _x\big|\pi_{n_0 \delta} 1_{(0, \infty)}(x) e^{-\tilde\alpha n_0 \delta}-n_1 V(x)\big|<\varepsilon
		\end{align*}
		by Proposition \ref{prop4.1}. We write $W_{\delta, k}^1=e^{-\tilde{\alpha} k n \delta_0} X_{k n \delta_0}(\infty)$, $k \in \mathbb{N} \backslash \{0\}$. Using (\ref{4.2}) we have
		\begin{align*}
			W_{\delta, k+1}^1
			=\ &e^{-\tilde{\alpha}(k+1) n \delta_0} X_{(k+1) n \delta_0}(\infty)\\
			=\ & W_{\delta,k}^1\Big\{\frac{1}{X_{kn\delta_0}(\infty)}\sum_{i=1}^{X_{kn\delta_0}(\infty)}e^{-\tilde{\alpha}n_0 \delta}\big[X_{n_0 \delta}^{x_i}(\infty)-\pi_{n_0 \delta}1_{(0,\infty)}(x_i)\big]\Big\}\\
			&+W_{\delta,k}^1\Big\{\frac{1}{X_{kn\delta_0}(\infty)}\sum_{i=1}^{X_{kn\delta_0}(\infty)}\big[e^{-\tilde{\alpha}n_0 \delta}\pi_{n_0 \delta}1_{(0,\infty)}(x_i)-n_1 V(x_i)\big]\Big\}\\
			&+n_1 e^{-\tilde{\alpha}n_0 \delta} \langle X_{kn_0 \delta},V\rangle.
		\end{align*}
		Then it follows from Proposition \ref{prop5.2'} that there exists a finite integer valued random variable $K_0$ such that for any $k\ge K_0$,
		\begin{align*}
			W_{\delta, k+1}^1 \le 2 \varepsilon W_{\delta, k}^1 +C_V,\quad \mathbf{P}\text{-a.s.}. 
		\end{align*}
		Iterating this proves that $\limsup _{k \rightarrow \infty} W_{\delta, k}^1<\infty$, $\mathbf{P}$-a.s.. By similar calculations we have
		\begin{align*}
			\limsup _{k \rightarrow \infty} e^{-\tilde\alpha\left(k n_0+j\right) \delta} X_{\left(k n_0+j\right) \delta}(\infty)<\infty,\quad \mathbf{P}\text {-a.s.},
		\end{align*}
		for $j=1,2, \cdots,(n_0-1)$. Notice that for any $n \in \mathbb{N} \backslash \{0\}$, there exists $j_0 \in\{1,2, \cdots,(n_0-1)\}$ such that $n/n_0=\lfloor n/n_0\rfloor+j_0/n_0$, where $\lfloor x\rfloor$ denotes the greatest integer $\le x$. Then we obtain
		\begin{align*}
			e^{-\tilde\alpha n\delta} X_{n\delta}(\infty)
			&=\exp \big\{-\tilde{\alpha}(\lfloor n / n_0\rfloor n_0+j_0) \delta\big\} X_{(\lfloor n / n_0\rfloor n_0+j_0) \delta}(\infty)\\
			&\le \sum_{j=0}^{n_0-1}\exp \big\{-\tilde{\alpha}(\lfloor n / n_0\rfloor n_0+j) \delta\big\} X_{(\lfloor n / n_0\rfloor n_0+j) \delta}(\infty),
		\end{align*}
		which follows the desired result.
		$\hfill\square$ 
	\end{prof}
	
	\begin{theorem}\label{th5.2'}
		Suppose that (\ref{2.*}), (\ref{super}), (\ref{3.3}) and (\ref{5.*}) hold. Then for any fixed  $\delta>0$ and $f\in B(0,\infty)^+$, conditioning on the non-extinction event, we have
		\begin{align*}
			A_{n\delta}(f)=\frac{\langle X_{n\delta},f\rangle}{X_{n\delta}(\infty)}\xrightarrow{\mathbf{P}\text{-a.s.}}A(f)<\infty,\quad \text{as}\ n\rightarrow\infty,
		\end{align*}
		where $A(f)$ is given by (\ref{4.*a}).
	\end{theorem}
	
	\begin{prof}
		Using (\ref{4.4}), Propositions \ref{prop4.1}, \ref{prop5.2'},  \ref{prop5.3} and the monotone convergence theorem we easily obtain the desired result.
		$\hfill\square$ 
	\end{prof}
	
	We next push the almost sure convergence on the lattice to the whole continuum.  Further, as an application we give the almost sure convergence of $W_t^{f}$.
	
	\begin{theorem}\label{th5.1''}
		Suppose that (\ref{2.*}), (\ref{super}), (\ref{3.3}) and (\ref{5.*}) hold. Then for any $x\in (0,\infty)$, conditioning on the non-extinction event, we have
		\begin{align*}
			A_{t}(x)\xrightarrow{\mathbf{P}\text{-a.s.}}A(x),\quad \text{as}\ t\rightarrow\infty,
		\end{align*}
		where $A(x)$ is given by (\ref{4.a}).
	\end{theorem}
	
	\begin{prof}
		For any $x\in(0,\infty)$, $1<\delta<x$ and $n\delta\le t<(n+1)\delta$ we have
		\begin{align}\label{5.8}
			X_{n \delta}(x+\delta)-\sum_{i=1}^{X_{n\delta}(\infty)} X_0^{x_i}(\delta) \le X_t(x) \le X_{(n+1) \delta}(x-\delta)+\sum_{i=1}^{X_{n\delta}(\infty)} X_0^{x_i}(\delta),
		\end{align}
		where $\{x_i; i=1,2,\ldots,X_{n\delta}(\infty)\}$ is the remaining-lifetime chart at time $n\delta$, $\big\{X_0^{x_i}; i=1,2,\ldots,X_t(\infty)\big\}$ are independent and further if $x_i=y$ then the conditional distribution of $X_0^{x_i}$ is the same as $X_0^{y}$ defined as (\ref{4.y}). It is easy to check that
		\begin{align}\label{5.9}
			A_t(x)\le\frac{\frac{X_{(n+1) \delta}(x-\delta)}{X_{(n+1) \delta}(\infty)} \cdot \frac{X_{(n+1) \delta}(\infty)}{X_{n\delta}(\infty)}+\frac{1}{X_{n\delta}(\infty)} \sum_{i=1}^{X_{n\delta}(\infty)} X_0^{x_i}(\delta)}{1-\frac{1}{X_{n\delta}(\infty)} \sum_{i=1}^{X_{n\delta}(\infty)} X_0^{x_i}(\delta)}.
		\end{align}
		By (\ref{4.2}) and Proposition \ref{prop5.2'}, as $n\rightarrow\infty$ we have
			\begin{align*}
				\frac{1}{X_{n \delta}(\infty)}&\Big\{X_{(n+1) \delta}(\infty)-\sum_{i=1}^{X_{n\delta}(\infty)} \pi_\delta 1_{(0, \infty)}(x_i)\Big\}\\
				&=\frac{1}{X_{n \delta}(\infty)}\sum_{i=1}^{X_{n\delta}(\infty)}\big\{X_{\delta}^{x_i}(\infty)- \pi_\delta 1_{(0, \infty)}(x_i)\big\}\xrightarrow{\mathbf{P}\text{-a.s.}}0
			\end{align*}
		and
		\begin{align*}
			\frac{1}{X_{n \delta}(\infty)}\sum_{i=1}^{X_{n\delta}(\infty)}\big\{X_{0}^{x_i}(\delta)- \pi_0 1_{(0, \delta)}(x_i)\big\}\xrightarrow{\mathbf{P}\text{-a.s.}}0.
		\end{align*}
		By Theorem \ref{th5.2'} we have
			\begin{align*}
				\lim_{n \rightarrow \infty}\frac{1}{X_{n \delta}(\infty)}\!\!\!\sum_{i=1}^{X_{n\delta}(\infty)}\!\!\!\pi_\delta 1_{(0, \infty)}(x_i)=\!\lim_{n \rightarrow \infty}\! A_{n\delta}(\pi_\delta 1_{(0, \infty)})=A(\pi_\delta 1_{(0, \infty)})=1+r_1(\delta),\quad \mathbf{P}\text{-a.s.},
			\end{align*}
		where
		\begin{align*}
			r_1(\delta)=\frac{\int_0^{\infty} e^{-\tilde\alpha u} d u \int_u^{\infty}\big[\pi_\delta 1_{(0, \infty)}(y-u)-1_{(0, \infty)}(y-u)\big] G(d y)}{\int_0^{\infty} e^{-\tilde\alpha u}[1-G(u)] d u}\rightarrow 0,\quad \text{as}\ \delta\downarrow 0.
		\end{align*}
		Notice that
		\begin{align*}
			\lim_{n \rightarrow \infty}\frac{1}{X_{n \delta}(\infty)}\!\!\sum_{i=1}^{X_{n\delta}(\infty)}\!\!\pi_0 1_{(0, \delta)}(x_i)=\lim_{n \rightarrow \infty} A_{n\delta}(\pi_0 1_{(0, \delta)})=A(\pi_0 1_{(0, \delta)})=A(\delta),\quad \mathbf{P}\text{-a.s.},
		\end{align*}
		where $A(\delta)$ is given by (\ref{4.a}). Since $A(\delta)\downarrow 0$ as $\delta\downarrow 0$. There exists $\delta_0>0$ such that $A(\delta)<1$ for any $\delta\in(0,\delta_0)$. Therefore, it follows from (\ref{5.9}) that for any $\delta\in(0,\delta_0\wedge x)$,
		\begin{align*}
			\limsup_{t \rightarrow \infty} A_t(x) \le \frac{A(x-\delta)\left(1+r_1(\delta)\right)+A(\delta)}{1-A(\delta)}<\infty.
		\end{align*}
		Letting $\delta\downarrow 0$ yields $\limsup_{t \rightarrow \infty} A_t(x) \le A(x)$. On the other hand, it is also easy to see that
		\begin{align*}
			A_t(x)\ge\frac{\frac{X_{n \delta}(x+\delta)}{X_{n\delta}(\infty)}-\frac{1}{X_{n\delta}(\infty)} \sum_{i=1}^{X_{n\delta}(\infty)} X_0^{x_i}(\delta)}{\frac{X_{(n+1) \delta}(\infty)}{X_{n\delta}(\infty)}+\frac{1}{X_{n\delta}(\infty)} \sum_{i=1}^{X_{n\delta}(\infty)} X_0^{x_i}(\delta)}.
		\end{align*}
		Then $\liminf_{t \rightarrow \infty} A_t(x) \ge A(x)$ can be obtained in a similar ways. Thus we get 
		the desired result.
		$\hfill\square$ 
	\end{prof}
	
	\begin{theorem}\label{th5.2''}
		Suppose that (\ref{2.*}), (\ref{super}), (\ref{3.3}) and (\ref{5.*}) hold. Then for any $f\in B(0,\infty)^+$, conditioning on the non-extinction event, we have
		\begin{align*}
			A_{t}(f)\xrightarrow{\mathbf{P}\text{-a.s.}}A(f),\quad \text{as}\ t\rightarrow\infty,
		\end{align*}
		where $A(f)$ is given by (\ref{4.*a}).
	\end{theorem}
	
	\begin{theorem}\label{th5.A}
		Suppose that (\ref{2.*}), (\ref{super}), (\ref{3.3}) and (\ref{5.*}) hold. Let $V\in B(0,\infty)^+$ be given by (\ref{4.v}). Then for any $f\in B(0,\infty)^+$, as $t\rightarrow \infty$, 
		$$W_t^f \xrightarrow{\mathbf{P}\text{-a.s.}} \widetilde{W}_{\infty}^{f}:=\frac{A(f)}{A(V)}\widetilde{W}_{\infty}^{V}$$ exists and is not identically zero, where
		$A(\cdot)$ is given as (\ref{4.*a}), $\widetilde{W}_{\infty}^{V}$ is given by (\ref{5.v}). Furthermore, the distribution of $\widetilde{W}_{\infty}^{f}$ is the same as $W_{\infty}^{f}$ given in Theorem \ref{th3.A}.
	\end{theorem}
	
	Naturally, if $\widetilde{W}_{\infty}^{f}$ exists and is not identically zero, then the absolute continuity of $\widetilde{W}_{\infty}^f$ is equivalent to this of $W_{\infty}^{f}$ given in Theorem \ref{th3.A}.

	\section{A central limit theorem}\label{section6}
	
	In this section we want to give a central limit theorem of $(X_t:t\ge 0)$. That is, for any $f\in B(0,\infty)^+$ we consider the convergence in distribution of 
	\begin{align*}
		\frac{\langle X_t, f\rangle-\mathbf{E}[\langle X_t, f\rangle]}{\sqrt{X_t(\infty)}}
	\end{align*}
	as $t\rightarrow\infty$. Naturally we just establish the results conditioning on the non-extinction event.
	
	It follows by arguments similar to (\ref{4.2}) that
	\begin{align}\label{6.e}
		\langle X_{t+s},f\rangle=\sum_{i=1}^{X_t(\infty)} \langle X_s^{x_i},f\rangle,\quad t,s\ge 0,
	\end{align}
	where $\{x_i; i=1,2,\ldots,X_t(\infty)\}$ is the remaining-lifetime chart at time $t$, $\big\{X_s^{x_i}; i=1,2,\ldots,X_t(\infty)\big\}$ are independent and further if $x_i=y$ then the conditional distribution of $X_s^{x_i}$ is the same as $X_s^{y}$ defined as (\ref{4.y}). Then
	\begin{align*}
		\mathbf{E}[\langle X_{t+s} , f\rangle]=\mathbf{E}\Big[\mathbf{E}\Big[\sum_{i=1}^{X_t(\infty)}\langle X_s^{x_i}, f\rangle\Big| \mathcal{F}_t\Big]\Big]\!=\mathbf{E}\Big[\sum_{i=1}^{X_t(\infty)}\!\!\mathbf{E}[\langle X_s^{x_i}, f\rangle]\Big]\!=\mathbf{E}\Big[\sum_{i=1}^{X_t(\infty)}\!\!\pi_{s}f(x_i)\Big].
	\end{align*}
	It follows that
	\begin{align*}
		\frac{\langle X_{t+s}, f\rangle-\mathbf{E}[\langle X_{t+s}, f\rangle]}{\sqrt{X_{t+s}(\infty)}}
		=\ &\sqrt{\frac{X_{t}(\infty)}{X_{t+s}(\infty)}}\frac{1}{\sqrt{X_{t}(\infty)}}\sum_{i=1}^{X_{t}(\infty)}\big[\langle X_s^{x_i},f\rangle-\pi_s f(x_i)\big]\\
		&+\frac{1}{\sqrt{X_{t+s}(\infty)}}\Big\{\sum_{i=1}^{X_{t}(\infty)}\pi_s f(x_i)-\mathbf{E}\Big[\sum_{i=1}^{X_{t}(\infty)}\pi_s f(x_i)\Big]\Big\}\\
		=\ &\sqrt{\frac{X_{t}(\infty)e^{-\tilde{\alpha}t}}{X_{t+s}(\infty)e^{-\tilde{\alpha}(t+s)}}}\frac{1}{\sqrt{X_{t}(\infty)}}\!\!\sum_{i=1}^{X_{t}(\infty)}\!\!\big[\langle X_s^{x_i},f\rangle-\pi_s f(x_i)\big]e^{-\frac{1}{2}\tilde{\alpha}s}\\
		&+\frac{1}{\sqrt{X_{t+s}(\infty)}}\Big\{\sum_{i=1}^{X_{t}(\infty)}\pi_s f(x_i)-\mathbf{E}\Big[\sum_{i=1}^{X_{t}(\infty)}\pi_s f(x_i)\Big]\Big\}\\
		=:\ &\sqrt{\frac{W_t^1}{W_{t+s}^1}}B_1(t,s)+B_2(t,s).
	\end{align*}
	Now it suffices to determine the asymptotic behavior of the random variables $B_1(t,s)$ and $B_2(t,s)$ as $t\rightarrow\infty$. To this end, we need to study the second-moment of $(X_t :t\ge 0)$. By Cheng and Li \cite[Proposition 2.5]{Cheng23+}, the proof of the following proposition is similar to that of Proposition \ref{prop2.1}.
	
	\begin{proposition}\label{prop6.g} 
		Suppose that $\|g''(\cdot,1-)\|<\infty$. Then for any $t \ge 0$ we have
		\begin{align}\label{6.1}
			\mathbf{E}[\langle X_{t}, f\rangle^2]=\langle G, \pi_{t} f\rangle^{2}+\langle G, \gamma_{t} f\rangle, \quad f \in B(0,\infty),
		\end{align}
		where $(\pi_{t})_{t \ge 0}$ is defined by (\ref{2.2'}) and $(t,x)\mapsto \gamma_{t}f(x)$ is the unique solution of
		\begin{align}\label{6.2}
			\gamma_{t} f(x)
			=&\int_{0}^{t} \alpha(x-s) \big[g''(x-s,1-) \langle G,\pi_{t-s}f\rangle^{2} + g'(x-s,1-) \langle G,(\pi_{t-s}f)^{2}\rangle\big] ds \nonumber\\
			&+\int_{0}^{t} \alpha(x-s) g'(x-s,1-) \langle G,\gamma_{t-s}f\rangle ds.
		\end{align}
	\end{proposition}
	
	\begin{proposition}\label{prop6.*}
		Suppose that (\ref{2.*}), (\ref{super}) hold and $\|g''(\cdot,1-)\|<\infty$. In addition, assume that the function 
			\begin{align}\label{pi}
				\Pi(t):=&\ e^{-\tilde{\alpha} t}\int_0^{\infty}G(dx)\int_{0}^{t}\alpha(x-r) g''(x-r,1-) \langle G,\pi_{t+s-r}f\rangle^{2}dr\nonumber\\
				&\ +e^{-\tilde{\alpha} t}\int_0^{\infty}G(dx)\int_{0}^{t}\alpha(x-r) g'(x-r,1-) \langle G,(\pi_{t+s-r}f)^{2}\rangle dr
			\end{align}
		is directly Riemann integrable over $[0,\infty)$. Then 
		\begin{align*}
			\lim _{s \rightarrow \infty} e^{-\tilde{\alpha} s} \lim _{t \rightarrow \infty} e^{-\tilde{\alpha}t}\langle G, \gamma_t(\pi_s f)\rangle=0.
		\end{align*}
	\end{proposition}
	
	\begin{prof}
		For any $t,s\ge 0$, $x\in (0,\infty)$ and $f\in B(0,\infty)^+$, We first observe that 
		\begin{align*}
			\pi_{t+s}f(x)=\hat{\mathbf{E}}_{x}[\langle \hat{X}_{t+s},f\rangle]=\hat{\mathbf{E}}_{x}\Big[\sum_{i=1}^{\hat{X}_{t}(\infty)}\!\!\hat{\mathbf{E}}_{x_i}[\langle \hat{X}_{t+s},f\rangle]\Big]\!=\hat{\mathbf{E}}_{x}[\langle \hat{X}_{t},\pi_{s}f\rangle]=\pi_{t}(\pi_s f)(x).
		\end{align*}
		Then by (\ref{6.2}) we have
		\begin{align*}
			e^{-\tilde{\alpha} t}\langle G, \gamma_t(\pi_s f)\rangle
			=\Pi(t)+\int_{0}^{t} e^{-\tilde{\alpha} (t-r)} \langle G,\gamma_{t-r}(\pi_s f)\rangle e^{-\tilde{\alpha} r} F(dr).
		\end{align*}
		Then by (\ref{2.*}) and the general result on defective renewal equation; see, e.g., Jagers \cite[Theorem 5.2.6]{Jagers75}, we have 
		\begin{align*}
			\lim _{t \rightarrow \infty} e^{-\tilde\alpha t}\langle G,\gamma_t(\pi_s f)\rangle
			&=\frac{1}{c_9}\int_0^{\infty}\Pi(u)du\\
			&=\frac{1}{c_9}\int_s^{\infty}e^{-\tilde{\alpha} (u-s)}du\int_0^{\infty}G(dx)\int_{0}^{u-s}\widetilde{\Pi}(u,x,r) dr,
		\end{align*}
		where $c_9=\int_0^{\infty} u e^{-\tilde\alpha u} F(d u)$ and
		\begin{align*}
			\widetilde{\Pi}(u,x,r)=\alpha(x-r) \big[g''(x-r,1-) \langle G,\pi_{u-r}f\rangle^{2} + g'(x-r,1-) \langle G,(\pi_{u-r}f)^{2}\rangle\big].
		\end{align*}
		Then as $s\rightarrow\infty$,
		\begin{align*}
			e^{-\tilde\alpha s}\lim _{t \rightarrow \infty} e^{-\tilde\alpha t}\langle G,\gamma_t(\pi_s f)\rangle
			&=\frac{1}{c_9}\int_s^{\infty}e^{-\tilde{\alpha} u}du\int_0^{\infty}G(dx)\int_{0}^{u-s}\widetilde{\Pi}(u,x,r) dr\\
			&\le\frac{1}{c_9}\int_s^{\infty}e^{-\tilde{\alpha} u}du\int_0^{\infty}G(dx)\int_{0}^{u}\widetilde{\Pi}(u,x,r) dr\\
			&\rightarrow 0,
		\end{align*}
		which implies the desired result.
		$\hfill\square$ 
	\end{prof}
	
	\begin{proposition}\label{prop6.1}
		Suppose that (\ref{2.*}), (\ref{super}) and (\ref{3.3}) hold. If $\|g''(\cdot,1-)\|<\infty$ and $\Pi(t)$ given as (\ref{pi}) is directly Riemann integrable over $[0,\infty)$, then for any $\varepsilon>0$ and $\delta>0$, there exists $s_0(\varepsilon,\delta)$ such that
		\begin{align*}
			\lim_{t\rightarrow \infty}\mathbf{P}[|B_2(t,s)|>\varepsilon]<\delta,\quad \text{for all}\ s\ge s_0(\varepsilon,\delta).
		\end{align*}
	\end{proposition}
	
	\begin{prof}
		By Theorem \ref{th5.A} we have $W_t^1 \xrightarrow{\mathbf{P}\text{-a.s.}} \widetilde{W}_{\infty}^{1}$ exists and conditioning on the non-extinction event, $\mathbf{P}[\widetilde{W}_{\infty}^{1}>0]=1$. Then for any $\delta>0$ there exists $x>0$ such that $\mathbf{P}[\widetilde{W}_{\infty}^{1}\le x]=\delta/3$ and for any $0<\varepsilon'<x/2$ there exists $s_0 '>0$ such that 
		\begin{align*}
			\mathbf{P}\big[|W_{t+s_0^{\prime}}^1-\widetilde{W}^1_{\infty}|>\varepsilon^{\prime}\big]<\frac{\delta}{3},\quad \text {for all}\ t \ge 0.
		\end{align*}
		Then for any $s\ge s_0'$ and $t\ge 0$ we have
		\begin{align}\label{6.8}
			\mathbf{P}[|B_2(t,s)|>\varepsilon]
			\le\ &\mathbf{P}\big[|B_2(t, x)|>\varepsilon, |W_{t+s}^1-\widetilde{W}_{\infty}^1|<\varepsilon^{\prime}, \widetilde{W}_{\infty}^1>x\big]\nonumber\\
			&+\mathbf{P}\big[|W_{t+s}^1-\widetilde{W}^1_{\infty}|\ge\varepsilon^{\prime}\big]+\mathbf{P}\big[\widetilde{W}_{\infty}^1\le x\big]\nonumber\\
			\le\ &\mathbf{P}\big[|B_2(t, x)|>\varepsilon, |W_{t+s}^1-\widetilde{W}_{\infty}^1|<\varepsilon^{\prime}, \widetilde{W}_{\infty}^1>x\big]+\frac{2\delta}{3}.
		\end{align}
		By Chebyshev's inequality,
		\begin{align}\label{6.9}
			\mathbf{P}&\big[|B_2(t, x)|>\varepsilon, |W_{t+s}^1-\widetilde{W}_{\infty}^1|<\varepsilon^{\prime}, \widetilde{W}_{\infty}^1>x\big]\nonumber\\
			&\le\mathbf{P}\Big[\Big|\frac{\langle X_t, \pi_s f\rangle-\mathbf{E}[\langle X_t, \pi_s f\rangle]}{\sqrt{X_{t+s}(\infty)}}\Big|>\varepsilon, X_{t+s}(\infty)>(x-\varepsilon')e^{\tilde{\alpha}(t+s)}\Big]\nonumber\\
			&\le\mathbf{P}\Big[|\langle X_t, \pi_s f\rangle-\mathbf{E}[\langle X_t, \pi_s f\rangle]|>\varepsilon\sqrt{x-\varepsilon'}e^{\frac{1}{2}\tilde{\alpha}(t+s)}\Big]\nonumber\\
			&\le\frac{e^{-\tilde{\alpha}(t+s)}}{\varepsilon^2(x-\varepsilon')} \operatorname{Var}[\langle X_t, \pi_s f\rangle]\nonumber\\
			&\le\frac{e^{-\tilde{\alpha}(t+s)}}{\varepsilon^2(x-\varepsilon')} \langle G,\gamma_{t}(\pi_{s} f)\rangle.
		\end{align}
		From Proposition \ref{prop6.*} it follows that there exists $s_0''$ such that 
		\begin{align}\label{6.10}
			e^{-\tilde\alpha s}\lim _{t \rightarrow \infty} e^{-\tilde\alpha t}\langle G,\gamma_t(\pi_s f)\rangle\le \frac{\delta}{3}(x-\varepsilon')\varepsilon^2,\quad \text{for all}\ t\ge 0,\ s\ge s_0''.
		\end{align}
		Let $s_0(\varepsilon,\delta)=\max\{s_0',s_0''\}$. In view of (\ref{6.8})-(\ref{6.10}), it is easy to see that
		\begin{align*}
			\lim_{t\rightarrow \infty}\mathbf{P}[|B_2(t,s)|>\varepsilon]<\delta,\quad \text{for all}\ s\ge s_0(\varepsilon,\delta),
		\end{align*}
		which completes the proof.
		$\hfill\square$ 
	\end{prof}
	
	\begin{proposition}\label{prop6.2}
		Suppose that (\ref{2.*}), (\ref{super}) and (\ref{3.3}) hold. If $\|g''(\cdot,1-)\|<\infty$, then for any fixed $s_0>0$, as $t\rightarrow\infty$,
		\begin{align*}
			\operatorname{Var}\big[B_1(t,s_0)\big|\mathcal{F}_t\big]\xrightarrow{\mathbf{P}\text{-a.s.}}e^{-\tilde{\alpha}s_0}A(\gamma_{s_0}f),
		\end{align*}
		where $A(\cdot)$ is given by  (\ref{4.*a}) and $(t,x)\mapsto \gamma_{t}f(x)$ is the unique solution of (\ref{6.2}).
	\end{proposition}
	
	\begin{prof}
		Write $Y_t^{x_i}(s_0)=[\langle X_{s_0}^{x_i}, f\rangle-\pi_{s_0} f(x_i)] e^{-\frac{1}{2} \tilde{\alpha} s_0}$, $i=1,2,\ldots,X_t(\infty)$. Then
		\begin{align*}
			B_1(t, s_0)=\frac{1}{\sqrt{X_t(\infty)}} \sum_{i=1}^{X_t(\infty)} Y_t^{x_i}(s_0).
		\end{align*}
		Since $\{Y_t^{x_i}(s_0); i=1,2,\ldots, X_t(\infty)\}$ are mutually independent conditioned on $\mathcal{F}_t$ and also independent of $X_t(\infty)$. Then it is easy to check that
		\begin{align*}
			\operatorname{Var}\big[B_1(t,s_0)\big|\mathcal{F}_t\big]=\frac{1}{\sqrt{X_t(\infty)}} \sum_{i=1}^{X_t(\infty)} \operatorname{Var}\big[Y_t^{x_i}(s_0)\big|\mathcal{F}_t^X\big].
		\end{align*}
		Notice that $\mathbf{E}\big[Y_t^{x_i}(s_0) \big| \mathcal{F}_t\big]=0$ and
		\begin{align*}
			\operatorname{Var}\big[Y_t^{x_i}(s_0)\big|\mathcal{F}_t\big]
			&=\mathbf{E}\big[(Y_t^{x_i}(s_0))^2 \big| \mathcal{F}_t\big]\\
			&=e^{-\tilde{\alpha} s_0} \mathbf{E}\big[\langle X_{s_0}^{x_i}, f\rangle^2 \big| \mathcal{F}_t\big]-e^{-\tilde{\alpha} s_0}[\pi_{s_0}f(x_i)]^2\\
			&=e^{-\tilde{\alpha} s_0}\gamma_{s_0}f(x_i).
		\end{align*}
		For any fixed $t\ge 0$, in view of (\ref{6.2}) we can use Gronwall's inequality to see $\gamma_{t}f(\cdot)\in B(0,\infty)^+$. Then by Theorem \ref{th5.2''} we get
		\begin{align*}
			\operatorname{Var}\big[B_1(t,s_0)\big|\mathcal{F}_t\big]
			&=\frac{1}{\sqrt{X_t(\infty)}} \!\sum_{i=1}^{X_t(\infty)} \!\!\operatorname{Var}\big[Y_t^{x_i}(s_0)\big|\mathcal{F}_t\big]
			=\frac{1}{\sqrt{X_t(\infty)}} \!\sum_{i=1}^{X_t(\infty)}\!\!e^{-\tilde{\alpha} s_0}\gamma_{s_0}f(x_i)\\
			&=e^{-\tilde{\alpha}s_0}A_t(\gamma_{s_0}f)
			\xrightarrow{\mathbf{P}\text{-a.s.}}e^{-\tilde{\alpha}s_0}A(\gamma_{s_0}f),
		\end{align*}
		as $t\rightarrow\infty$.
		$\hfill\square$ 
	\end{prof}
	
	\begin{proposition}\label{prop6.3}
		Suppose that $\|g''(\cdot,1-)\|<\infty$. Then for any fixed $s_0 \ge 0$ and $\delta>0$, as $t\rightarrow\infty$ we have
		\begin{align*}
			\sup _{y\ge 0} \mathbf{E}\Big\{[Y_t^y(s_0)]^2 ;|Y_t^y(s_0)|>\delta e^{\frac{1}{2} \tilde\alpha t}\Big\} \rightarrow 0,
		\end{align*}
		where $Y_t^{y}(s_0)=[\langle X_{s_0}^{y}, f\rangle-\pi_{s_0} f(y)] e^{-\frac{1}{2} \tilde{\alpha} s_0}$ given as in the proof of Proposition \ref{prop6.2}.
	\end{proposition}
	
	\begin{prof}
		Since $\sup _{y\ge0} \pi_{s_0} f(y)<\infty$ and
		\begin{align*}
			\mathbf{E}\Big\{[Y_t^y(s_0)]^2 ;|Y_t^y(s_0)|>\delta e^{\frac{1}{2} \tilde\alpha t}\Big\}\le\ 
			&2 \mathbf{E}\Big\{e^{-\tilde{\alpha} s_0}\langle X_{s_0}^y, f\rangle^2 ; \langle X_{s_0}^y, f\rangle>\delta e^{\frac{1}{2} \tilde{\alpha}(t+s_0)}\Big\}\\
			&+2 \mathbf{E}\Big\{e^{-\tilde{\alpha} s_0}(\pi_{s_0}f(y))^2 ; \langle X_{s_0}^y, f\rangle>\delta e^{\frac{1}{2} \tilde{\alpha}(t+s_0)}\Big\}.
		\end{align*}
		Now it suffices to show that as $t\rightarrow\infty$,
		\begin{align}\label{6.i}
			\sup_{y\ge 0}\mathbf{P}\Big[\langle X_{s_0}^y, f\rangle>\delta e^{\frac{1}{2} \tilde{\alpha}(t+s_0)}\Big]\rightarrow 0
		\end{align}
		and 
		\begin{align}\label{6.ii}
			\sup_{y\ge 0}e^{-\tilde{\alpha} s_0}\mathbf{E}\Big\{\langle X_{s_0}^y, f\rangle^2 ; \langle X_{s_0}^y, f\rangle>\delta e^{\frac{1}{2} \tilde{\alpha}(t+s_0)}\Big\}\rightarrow 0.
		\end{align}
		Indeed, 
		\begin{align*}
			\sup_{y\ge 0}\mathbf{P}\Big[\langle X_{s_0}^y, f\rangle>\delta e^{\frac{1}{2} \tilde{\alpha}(t+s_0)}\Big]
			&\le \delta^{-1} e^{-\frac{1}{2} \tilde{\alpha}(t+s_0)} \sup_{y\ge 0}\mathbf{E}\big[\langle X_{s_0}^y, f\rangle\big]\\
			&= \delta^{-1} e^{-\frac{1}{2} \tilde{\alpha}(t+s_0)} \sup_{y\ge 0}\|\pi_{s_0}f\|\\
			&\rightarrow 0,\quad\text{as}\ t\rightarrow\infty,
		\end{align*}
		which proves (\ref{6.i}). Turning to (\ref{6.ii}), by arguments similar to (\ref{5.0}) that $\langle X_{s_0}^y,f\rangle$ can be represented as the sum 
		\begin{align*}
			\langle X_{s_0}^y,f\rangle=\sum_{j=1}^{\eta_{\tau_1(y)}} \langle X_{s_0-\tau_1(y)}^{(j)},f\rangle +f(y-s_0),
		\end{align*}
		where $(X_{t}^{(j)}: t\ge 0)$,  $j=1,2,\ldots,\eta_{\tau_1(y)}$ are independent and have the same distribution as $(X_{t}: t\ge 0)$. Write $\Phi:=\sum_{j=1}^{\eta_{\tau_1(y)}} \langle X_{s_0-\tau_1(y)}^{(j)},f\rangle$. Then for any $\varepsilon>0$ we have
		\begin{align}\label{6.14}
			\mathbf{P}\big[\langle X_{s_0}^y, f\rangle^2 \ge \varepsilon\big] \le \mathbf{P}[\Phi^2 \ge \varepsilon]+\mathbf{P}\big[(f(y-s_0))^2 \ge \varepsilon\big].
		\end{align}
		It follows from (\ref{6.i}) that $\text{as}\ t\rightarrow\infty$,
		\begin{align}\label{6.15}
			\sup _{y \ge 0} \mathbf{E}\Big\{(f(y-s_0))^2 ;\langle X_{s_0}^y, f\rangle>\delta e^{\frac{1}{2} \tilde{\alpha}(t+s_0)}\Big\}
			&\le\|f\|^2\sup_{y\ge 0}\mathbf{P}\Big[\langle X_{s_0}^y, f\rangle>\delta e^{\frac{1}{2} \tilde{\alpha}(t+s_0)}\Big]\nonumber\\
			&\rightarrow 0.
		\end{align}
		On the other hand, notice that $\mathbf{P}[\tau_1(y) \in d s]=\alpha(y-s) e^{-\int_0^s \alpha(y-r) d r} d s$ and $\mathbf{P}[\eta_s=n]=p(y-s, n)$ for $s>0$ and $n\in\mathbb{N}$. Then
		\begin{align*}
			\mathbf{E}\Phi^2
			&=\int_0^{\infty} \sum_{n \in \mathbb{N}} b_8(y,s,n) \mathbf{E}\Big[\sum_{j=1}^n \langle X_{s_0-s}^{(j)}, f\rangle\Big]^2 d s\\
			&=\int_0^{\infty} \sum_{n \in \mathbb{N}} b_8(y,s,n) \Big\{\operatorname{Var}\Big[\sum_{j=1}^n \langle X_{s_0-s}^{(j)}, f\rangle\Big]+\Big[\sum_{j=1}^n \mathbf{E}\big[\langle X_{s_0-s}^{(j)}, f\rangle\big]\Big]^2\Big\} d s\\
			&=\int_0^{\infty} \sum_{n \in \mathbb{N}} b_8(y,s,n) \Big\{\sum_{j=1}^n \operatorname{Var}\big[\langle X_{s_0-s}^{(j)}, f\rangle\big]+\Big[\sum_{j=1}^n \mathbf{E}\big[\langle X_{s_0-s}^{(j)}, f\rangle\big]\Big]^2\Big\} d s\\
			&=n\langle G, \gamma_{s_0-s} f\rangle+n^2\langle G, \pi_{s_0-s} f\rangle^2 <\infty,
		\end{align*}
		where $b_8(y,s,n)=\alpha(y-s) p(y-s, n) e^{-\int_0^s \alpha(y-r) d r}$. Notice that
		\begin{align*}
			\sup _{y \ge 0} \mathbf{E}\big\{\Phi^2 ;\langle X_{s_0}^y, f\rangle>\delta e^{\frac{1}{2} \tilde{\alpha}(t+s_0)}\big\}
			\le\mathbf{E}\Big\{\Phi^2\Big[1_{\big\{\!\Phi>\delta e^{\frac{1}{2} \tilde{\alpha}(t+s_0)}\!\big\}}+1_{\big\{\!\|f\|>\delta e^{\frac{1}{2} \tilde{\alpha}(t+s_0)}\!\big\}}\Big]\Big\}.
		\end{align*}
		Since $\lim_{t\rightarrow \infty}\mathbf{P}[\Phi>\delta e^{\frac{1}{2} \tilde{\alpha}(t+s_0)}]=0$ and $\lim_{t\rightarrow \infty}\mathbf{P}[\|f\|>\delta e^{\frac{1}{2} \tilde{\alpha}(t+s_0)}]=0$ we conclude that as $t\rightarrow\infty$, 
		\begin{align}\label{6.16}
			\sup _{y \ge 0} \mathbf{E}\big\{\Phi^2 ;\langle X_{s_0}^y, f\rangle>\delta e^{\frac{1}{2} \tilde{\alpha}(t+s_0)}\big\}\rightarrow 0
		\end{align}
		by dominated convergence. Then (\ref{6.14})-(\ref{6.16}) imply that (\ref{6.ii}) holds.
		$\hfill\square$ 
	\end{prof}
	
	\begin{proposition}\label{prop6.4}
		Suppose that $\|g''(\cdot,1-)\|<\infty$. Then for any fixed $s_0 \ge 0$ and $\delta>0$, as $t\rightarrow\infty$ we have
		\begin{align*}
			\frac{1}{X_t(\infty)}\sum_{i=1}^{X_t(\infty)}\mathbf{E}\Big\{[Y_t^{x_i}(s_0)]^2 ;|Y_t^{x_i}(s_0)|>\delta \sqrt{X_t(\infty)}\Big|\mathcal{F}_t\Big\} \xrightarrow{\mathbf{P}} 0,
		\end{align*}
		where $Y_t^{x_i}(s_0)=[\langle X_{s_0}^{x_i}, f\rangle-\pi_{s_0} f(x_i)] e^{-\frac{1}{2} \tilde{\alpha} s_0}$ given as in the proof of Proposition \ref{prop6.2}.
	\end{proposition}
	
	\begin{prof}
		Write $S_t(s_0,\delta):=\frac{1}{X_t(\infty)}\sum_{i=1}^{X_t(\infty)}\mathbf{E}\big\{[Y_t^{x_i}(s_0)]^2 ;|Y_t^{x_i}(s_0)|>\delta \sqrt{X_t(\infty)}\big|\mathcal{F}_t\big\}$. Given $\delta_1>0$ and $\delta_2>0$, there exists $t_0>0$ and a set $\Omega_1\subset\Omega$ such that
		\begin{align}\label{6.17}
			\mathbf{P}(\Omega_1)>1-\delta_1
		\end{align}
		and
		\begin{align}\label{6.18}
			X_t(\infty)(\omega)\ge \delta_2 e^{\tilde{\alpha}t},\quad\text{for all}\ t>t_0, \omega\in\Omega_1. 
		\end{align}
		Then for any $\varepsilon>0$, by (\ref{6.17}) we have
		\begin{align*}
			\mathbf{P}[S_t(s_{0},\delta)>\varepsilon]
			&=\mathbf{P}[S_t(s_{0},\delta)>\varepsilon ; \Omega_1]+\mathbf{P}[S_t(s_{0},\delta)>\varepsilon ; \Omega\backslash \Omega_1]\\
			&\le \mathbf{P}[S_t(s_{0},\delta)>\varepsilon ; \Omega_1]+\delta_1.
		\end{align*}
		And for any $t\ge t_0$, by (\ref{6.18}) and Proposition \ref{prop6.3} we have
		\begin{align*}
			\mathbf{P}[S_t(s_{0},\delta)>\varepsilon ; \Omega_1]
			&\le\mathbf{P}\Big[\frac{1}{X_t(\infty)}\sum_{i=1}^{X_t(\infty)}\mathbf{E}\Big\{[Y_t^{x_i}(s_0)]^2 ;|Y_t^{x_i}(s_0)|>\delta \sqrt{\delta_2}e^{\frac{1}{2}\tilde{\alpha}t}\Big|\mathcal{F}_t\Big\}>\varepsilon\Big]\\
			&\le\mathbf{P}\Big[\sup_{x_i \ge 0}\mathbf{E}\Big\{[Y_t^{x_i}(s_0)]^2 ;|Y_t^{x_i}(s_0)|>\delta \sqrt{\delta_2}e^{\frac{1}{2}\tilde{\alpha}t}\Big|\mathcal{F}_t\Big\}>\varepsilon\Big]\\
			&\rightarrow 0,
		\end{align*}
		as $t\rightarrow\infty$. Then we conclude that $\lim_{t\rightarrow \infty}	\mathbf{P}[S_t(s_{0},\delta)>\varepsilon]<\delta_1$. Letting $\delta_1\downarrow 0$ we get the desired result.
		$\hfill\square$ 
	\end{prof}
	
	\begin{proposition}\label{prop6.5}
		Suppose that $\|g''(\cdot,1-)\|<\infty$. Let $A(\cdot)$ be given by (\ref{4.*a}) and $(t,x)\mapsto \gamma_{t}f(x)$ be the unique solution of (\ref{6.2}). Then for any fixed $s_0>0$, as $t\rightarrow\infty$,
		\begin{align*}
			B_1(t,s_0)\xrightarrow{d} \mathcal{N}(0,e^{-\tilde{\alpha}s}A(\gamma_{s}f)),
		\end{align*}
		where $\mathcal{N}(0,e^{-\tilde{\alpha}s}A(\gamma_{s}f))$ is the Gaussian distribution with mean zero and variance $e^{-\tilde{\alpha}s_0}A(\gamma_{s_0}f)$.
	\end{proposition}
	
	\begin{prof}
		For any $0<\theta<\infty$ we have
		\begin{align*}
			\mathbf{E}[\exp \{-\theta B_1(t, s_0)\} |\mathcal{F}_t]
			=\prod_{i=1}^{X_t(\infty)} \mathbf{E}\Big[\exp \Big\{-\theta \frac{Y_t^{x_i}(s_0)}{\sqrt{X_t(\infty)}}\Big\} \Big| \mathcal{F}_t\Big]=:\prod_{i=1}^{X_t(\infty)}\phi_{t}^{x_i}(s_0,\theta).
		\end{align*}
		As in the proof of the Lindeberg-Feller central limit theorem; see, e.g., Durrett \cite[Theorem 3.4.10]{Durrett19}, it is simple to show that as $t\rightarrow\infty$,
		\begin{align*}
			\prod_{i=1}^{X_t(\infty)}\phi_{t}^{x_i}(s_0,\theta)\xrightarrow{\mathbf{P}}\exp\Big\{\frac{\theta^2}{2}e^{-\tilde{\alpha}s_0}A(\gamma_{s_0}f)\Big\}
		\end{align*}
		by Propositions \ref{prop6.2} and \ref{prop6.4}. Notice that
		\begin{align*}
			\Big|\prod_{i=1}^{X_t(\infty)}\phi_{t}^{x_i}(s_0,\theta)-\exp\Big\{\frac{\theta^2}{2}e^{-\tilde{\alpha}s_0}A(\gamma_{s_0}f)\Big\}\Big|\le 2.
		\end{align*}
		Then we obtain that as $t\rightarrow\infty$,
			\begin{align*}
				\mathbf{E}&[\exp \{-\theta B_1(t, s_0)\}]\\
				&=\mathbf{E}\big\{\mathbf{E}[\exp \{-\theta B_1(t, s_0)\} |\mathcal{F}_t]\big\}\\
				&=\mathbf{E}\Big\{\mathbf{E}\Big[e^{-\theta B_1(t, s_0)}-\exp\Big\{\frac{\theta^2}{2}e^{-\tilde{\alpha}s_0}A(\gamma_{s_0}f)\Big\} \Big|\mathcal{F}_t\Big]\Big\}+\exp\Big\{\frac{\theta^2}{2}e^{-\tilde{\alpha}s_0}A(\gamma_{s_0}f)\Big\}\\
				&\rightarrow \exp\Big\{\frac{\theta^2}{2}e^{-\tilde{\alpha}s_0}A(\gamma_{s_0}f)\Big\},
			\end{align*}
		by dominated convergence, which completes the proof.
		$\hfill\square$ 
	\end{prof}
	
	\begin{proposition}\label{prop6.6}
		Suppose that  $\|g''(\cdot,1-)\|<\infty$ and $\Pi(t)$ given in (\ref{pi}) is directly Riemann integrable over $[0,\infty)$. Let $D_f:=\lim _{t \rightarrow \infty} e^{-\tilde\alpha t}\langle G,\gamma_t f\rangle$ and
		\begin{align*}
			\sigma(x)=\int_0^x \alpha(x-r) g'(x-r, 1-) e^{-\tilde{\alpha} r} d r,\quad x\in(0,\infty).
		\end{align*}
		Then we have
		\begin{align}\label{6.s}
			\lim_{s \rightarrow \infty}e^{-\tilde{\alpha}s}A(\gamma_{s}f)=A(\sigma)D_f,
		\end{align}
		where $A(\cdot)$ is given by  (\ref{4.*a}), $(t,x)\mapsto \gamma_{t}f(x)$ is the unique solution of (\ref{6.2}).
	\end{proposition}
	
	\begin{prof}
		Let $\widetilde{\Pi}(\cdot, \cdot, \cdot)$ be given as in the proof of Proposition \ref{prop6.*}. By (\ref{6.2}) we have
		\begin{align*}
			e^{-\tilde{\alpha} s}\langle G, \gamma_s f\rangle
			=e^{-\tilde{\alpha} s}\!\!\int_0^{\infty}\!\!G(dx)\!\int_{0}^{s}\!\! \widetilde{\Pi}(s, x, r) dr
			+\int_{0}^{s} e^{-\tilde{\alpha} (s-r)} \langle G,\gamma_{s-r}f\rangle e^{-\tilde{\alpha} r} F(dr).
		\end{align*}
		As in the proof of Proposition \ref{prop6.*} it is simple to see
		\begin{align*}
			D_f=\lim _{s \rightarrow \infty} e^{-\tilde\alpha s}\langle G,\gamma_s f\rangle
			=\frac{1}{c_9}\int_0^{\infty}e^{-\tilde{\alpha} u}du\int_0^{\infty}G(dx)\int_{0}^{u}\widetilde{\Pi}(u, x, r) dr<\infty,
		\end{align*}
		where $c_9$ and $\widetilde{\Pi}(u, x, r)$ are given as in the proof of Proposition \ref{prop6.*}. Then for any $s\ge 0$, there exists $b_7(f)\in(0,\infty)$ such that $e^{-\tilde\alpha s}\langle G,\gamma_s f\rangle\le b_7(f)$. Notice that for any $s\ge 0$,
		\begin{align*}
			\gamma_s f(y-u)
			=\int_{0}^{s}\! \widetilde{\Pi}(s, y-u, r) dr
			+\int_{0}^{s}\!\! \alpha(y-u-r) g'(y-u-r,1-) \langle G,\gamma_{s-r}f\rangle dr,
		\end{align*}
		and as $s\rightarrow\infty$,
		\begin{align*}
			\int_0^{\infty}\!\!\!e^{-\tilde{\alpha}u}du\!\int_u^{\infty}\!\!\!G(dy)e^{-\tilde{\alpha}s}\!\int_{0}^{s}\!\widetilde{\Pi}(s, y-u, r) dr
			&=\int_0^{\infty}\!\!\!e^{-\tilde{\alpha} (u+s)}du\!\int_0^{\infty}\!\!\!G(dy)\!\int_{u}^{u+s}\!\!\widetilde{\Pi}(u+s, y, r) dr\\
			&\le\int_0^{\infty}\!\!\!e^{-\tilde{\alpha} (u+s)}du\!\int_0^{\infty}\!\!\!G(dy)\!\int_{0}^{u+s}\!\!\widetilde{\Pi}(u+s, y, r) dr\\
			&=\int_s^{\infty}\!\!e^{-\tilde{\alpha} u}du\!\int_0^{\infty}\!\!G(dy)\!\int_{0}^{u}\!\!\widetilde{\Pi}(u, y, r) dr\\
			&\rightarrow 0.
		\end{align*}
		Then using dominated convergence we have
		\begin{align*}
			\lim_{s\rightarrow\infty}&e^{-\tilde{\alpha}s}A(\gamma_{s}f)\\
			&=\lim_{s\rightarrow\infty}\frac{\int_0^{\infty}\!\!e^{-\tilde{\alpha}u}du\!\int_u^{\infty}\!\!G(dy)\!\int_{0}^{s}\!\!\alpha(y-u-r) g'(y-u-r,1-) e^{-\tilde{\alpha}(s-r)}\langle G,\gamma_{s-r}f\rangle e^{-\tilde{\alpha}r}dr}{\int_0^{\infty}e^{-\tilde{\alpha}u}[1-G(u)]du}\\
			&= \frac{\int_0^{\infty}\!\!e^{-\tilde{\alpha}u}du\!\int_u^{\infty}\!\!G(dy)\!\int_{0}^{y-u}\!\!\alpha(y-u-r) g'(y-u-r,1-) e^{-\tilde{\alpha}r}dr}{\int_0^{\infty}e^{-\tilde{\alpha}u}[1-G(u)]du}\cdot D_f\\
			&=A(\sigma)D_f,
		\end{align*}
		where
		\begin{align*}
			\sigma(x)=\int_0^x \alpha(x-r) g'(x-r, 1-) e^{-\tilde{\alpha} r} d r.
		\end{align*}
		$\hfill\square$ 
	\end{prof}
	
	\begin{theorem}\label{th6.C}
		Suppose that (\ref{2.*}), (\ref{super}) and (\ref{3.3}) hold. If $\|g''(\cdot,1-)\|<\infty$ and $\Pi(t)$ given as (\ref{pi}) is directly Riemann integrable over $[0,\infty)$, then for any $f\in B(0,\infty)^+$, as $t\rightarrow\infty$,
		\begin{align*}
			\frac{\langle X_t, f\rangle-\mathbf{E}[\langle X_t, f\rangle]}{\sqrt{X_t(\infty)}}\xrightarrow{d} \mathcal{N}(0,A(\sigma)D_f),
		\end{align*} 		
		where $\mathcal{N}(0,A(\sigma)D_f)$ is the Gaussian distribution with mean zero and variance $A(\sigma)D_f$, $A(\cdot)$ is given by (\ref{4.*a}), $\sigma(\cdot)$ and $D_f$ are given in Proposition \ref{prop6.6}.
	\end{theorem}
	
	\begin{prof}
		Let $\mathcal{N}(0,b)$ be the Gaussian distribution with mean zero and variance $b$. Then letting $F_1(x)$ and $F_2(x)$ be the distribution of $\mathcal{N}(0,e^{-\tilde{\alpha}s}A(\gamma_{s}f))$ and $\mathcal{N}(0,A(\sigma)D_f)$, respectively. For any $\varepsilon>0$ and fixed $y\ge 0$, there exists $\delta_{\varepsilon}>0$ such that
		\begin{align}\label{6.20}
			|F_2(y+\delta_{\varepsilon})-F_2(y-\delta_{\varepsilon})|<\frac{\varepsilon}{3}.
		\end{align}
		Since $\lim_{s \rightarrow \infty}e^{-\tilde{\alpha}s}A(\gamma_{s}f)=A(\sigma)D_f$, there exists $s_1(\varepsilon)>0$ such that for any $s>s_1(\varepsilon)$,
		\begin{align}\label{6.21}
			|F_2(y+\delta_{\varepsilon})-F_1(y+\delta_{\varepsilon})|<\frac{\varepsilon}{3}\quad\text{and}\quad |F_2(y-\delta_{\varepsilon})-F_1(y-\delta_{\varepsilon})|<\frac{\varepsilon}{3}.
		\end{align}
		Let $s^*=\max\{s_0(\varepsilon,\delta_{\varepsilon}),s_1(\varepsilon)\}$, where $s_0(\varepsilon,\delta_{\varepsilon})$ is given in Proposition \ref{prop6.1}. Let $F_{1^*}(x)$ denotes the distribution of $\mathcal{N}(0,e^{-\tilde{\alpha}s^*}A(\gamma_{s^*}f))$. Then
		\begin{align}\label{6.22}
			&\limsup_{t \rightarrow \infty} \mathbf{P}\big[B_1(t, s^*)+B_2(t, s^*) \leq y\big]\nonumber\\
			&\qquad\le \limsup_{t \rightarrow \infty} \mathbf{P}\big[B_1(t, s^*)\le y+\delta_{\varepsilon}\big]+\limsup_{t \rightarrow \infty} \mathbf{P}\big[|B_2(t, s^*)|\ge  \delta_{\varepsilon}\big]\nonumber\\
			&\qquad\le F_{1^*}(y+\delta_{\varepsilon})+\frac{\varepsilon}{3}
		\end{align}
		and
		\begin{align}\label{6.23}
			&\liminf_{t \rightarrow \infty} \mathbf{P}\big[B_1(t, s^*)+B_2(t, s^*) \leq y\big]\nonumber\\
			&\qquad\ge \liminf_{t \rightarrow \infty} \mathbf{P}\big[B_1(t, s^*)\le y-\delta_{\varepsilon}\big]+\liminf_{t \rightarrow \infty} \mathbf{P}\big[|B_2(t, s^*)|\ge  \delta_{\varepsilon}\big]\nonumber\\
			&\qquad\ge F_{1^*}(y-\delta_{\varepsilon})-\frac{\varepsilon}{3}.
		\end{align}
		Notice that as $t\rightarrow\infty$,
		\begin{align*}
			\sqrt{W_t^1 / W_{t+s^*}^1} \xrightarrow{\mathbf{P}\text{-a.s.}} 1.
		\end{align*}
		Using this and (\ref{6.20})-(\ref{6.23}) we have
		\begin{align*}
			F_{2}(y)-\varepsilon\le\lim_{t\rightarrow \infty}\Big[\frac{\langle X_t, f\rangle-\mathbf{E}[\langle X_t, f\rangle]}{\sqrt{X_t(\infty)}}\le y\Big]\le F_2 (y)+\varepsilon.
		\end{align*}
		Since $\varepsilon>0$ is arbitrary, the proof is completed.
		$\hfill\square$ 
	\end{prof}

\end{document}